\documentclass[12pt]{article}

\usepackage{booktabs} 

\usepackage{color}

\usepackage{ifthen}
\usepackage{comment} 

\excludecomment{versionC}

\usepackage{amsmath}
\usepackage{amssymb}
\usepackage{amsfonts}
\usepackage{amsthm}
\usepackage{latexsym}
\usepackage{epsfig}
\usepackage{array}
\usepackage{boxedminipage}

\makeatletter
\renewcommand{\section}{\@startsection
{section} {1} {0mm} {-\baselineskip} {0.5\baselineskip}
{\large\bf}}
\renewcommand{\subsection}{\@startsection
{subsection} {2} {0mm} {-\baselineskip} {0.5\baselineskip}
{\normalsize\bf}} \makeatother

\newtheorem{defin}{Definition}
\newtheorem{theor}{Theorem}
\newtheorem{lema}{Lemma}
\newtheorem{propo}{Proposition}

\newtheorem{coro}{Corollary}

\newcommand{\vs}{\vspace{1mm}}
\newcommand{\vv}{\vspace{2mm}}

\newcommand{\vvvvv}{\vspace{5mm}}

\newcommand{\dsp}{\displaystyle}
\newcommand{\Csp}{Core-semiperiphery-periphery }
\newcommand{\csp}{core-semiperiphery-periphery }
\newcommand{\Cp}{Core-periphery }
\newcommand{\cp}{core-periphery }
\newcommand{\spe}{semiperiphery-periphery }
\newcommand{\s}{semiperiphery }
\newcommand{\p}{periphery }

\begin{document}

\begin{center}
{\Large\bf Twin subgraphs and \vspace{2mm} \\
core-semiperiphery-periphery structures\footnote{This is the
author's version of a paper accepted for publication in 
{\em Complexity}, 2018 (in press).}
}

\vs

{\sc 
Ricardo Riaza\footnote{Supported by Research Project 
MTM2015-67396-P (MINECO/FEDER). Email: {\em ricardo.riaza@upm.es}}}\\ 
Depto.\ de 
Matem\'{a}tica Aplicada 
a las TIC \&  \\
Information Processing and Telecommunications Center \\ 
ETSI 
Telecomunicaci\'{o}n, Universidad Polit\'{e}cnica de Madrid, Spain \\ 
\end{center}

\begin{abstract}
A standard approach to reduce the complexity of very large networks is 
to group together sets of nodes into clusters
according to some criterion which reflects certain structural properties
of the network.
Beyond the well-known modularity measures defining
communities, there are criteria based on the existence of similar
or identical connection patterns 
of a node or sets of nodes to the remainder of 
the network; this approach supports so-called positional analyses and the definition
of certain structures in social, commercial and economic networks.
A key notion in this context is that of {\em structurally equivalent}
or {\em twin} nodes, displaying exactly the same connection pattern to 
the remainder of the network. 

The first goal of this paper is to extend
this idea to subgraphs of arbitrary order of a given network,
by means of the notions of T-twin and F-twin 
subgraphs.
This research, which leads to graph-theoretic results
of independent interest, is motivated by the need to provide
a systematic approach to the analysis of \csp (CSP) structures, a notion
which is widely used in network theory but that somehow lacks a formal
treatment in the literature. 
The goal is to provide an analytical framework accommodating and
extending the idea that the unique (ideal) core-periphery (CP) structure 
is a 2-partitioned $K_2$, a fact which is here understood to rely
on the true-twin and false-twin notions for 
vertices already known in network theory. 
We provide a formal definition of such CSP structures in terms
of core eccentricities and periphery degrees, with semiperiphery
vertices acting as intermediaries between both. The 
T-twin and F-twin 
notions then make it possible to reduce the large number of resulting
structures by identifying isomorphic substructures which share
the connection pattern to the remainder of the graph, paving the
way for the decomposition and enumeration of CSP structures.
We compute explicitly the resulting CSP structures up 
to order six. 

We illustrate the scope of our results
by analyzing
a subnetwork of the well-known 
network of 
metal manufactures trade
arising
from 1994 world trade statistics.
As this example suggests, our approach 
can be naturally applied in complex
network theory and
seem to have many potential extensions, since 
the analytical properties of twin subgraphs and 
the
structure of 
CSP and other partitioned 
graphs admit further study.


\end{abstract}

\vspace{1mm}


\vvvvv \newpage

\section{Introduction}
\label{sec-intr}

The notion of a \cp (CP) structure can be traced back at least to some research 
on 
economic and commercial networks developed in the late 1970s and early 1980s
\cite{chasedunn75, nemeth85, snyder79}, largely emanating from the
influential work of 
Wallerstein on world systems analysis \cite{wallerstein74}. These ideas were revisited
and addressed in a more formal framework
 by Borgatti and Everett 
in \cite{borgatti99}. For these authors,
the two key ideas in the definition of a \cp structure 
in a network context
are those of 
a dense, cohesive core of heavily interconnected nodes
and a sparse periphery of nodes, essentially lacking any connections among them;
by contrast, the connection pattern
between the core and the periphery admits several definitions
and, actually, the core-periphery connection densities 
differ from some models to others. 
In idealized models, core nodes are fully connected among them, periphery nodes
are isolated (within the periphery subnetwork), whereas 
the core and the periphery may
either be fully connected or totally disconnected.
Since then, a great deal of research has been directed to
the detection of such core-periphery structures in real networks, measuring how
well they approximate the ideal ones, and to the
development of analytical and computational 
tools to classify nodes in such networks
(cf.\ \cite{borgattiBook, daSilva, dellarossa, hidalgo,
holme05, rombach14, zhangNewman} and references
therein).
Other approaches to the definition of a 
core-periphery structure can be found in \cite{chartrand00, csermely, gamble}. 

Even though the idea of a \csp (CSP) structure can
be also found in the 
aforementioned sociological 
works (cf.\ \cite{snyder79, wallerstein74}), and despite the fact
that this concept has been 
widely used since then 
(see e.g.\ \cite{dellarossa, deNooy, lazega, rombach14}),
the network literature seems to lack
a formal definition and a systematic classification of these CSP structures.
In the aforementioned paper by Borgatti and Everett \cite{borgatti99},
these authors indicate that there are many reasonable options 
to define a CSP structure and, further, discrete partitions with 
more than three classes.
The difficulty does not seem to rely on providing a formal definition
but on classifying the resulting ``reasonable options'', 
quoting these authors;
more precisely, there is a need for a notion of similar 
or equivalent subgraphs making it possible to somehow reduce the
number of different CSP structures.
When dealing with \cp structures, there
is a well-known subgraph similarity notion which makes this
reduction feasible, namely, that of {\em structural equivalence}
defining so-called {\em twin nodes} 
(broadly, two vertices are twins if they have the same neighbors;
a distinction is made between {\em true twins} and {\em false twins}
depending on whether both vertices are adjacent or not; details are given in Section \ref{sec-bac}).
Essentially, under structural equivalence, $K_2$ will be the unique 
core-periphery structure: 
details are provided later, but the
reader can think for the moment e.g.\ in the star $S_n$ as a network
with a unique core (the central node) to whom $n-1$ peripheric nodes
are attached;
all $n-1$ leaves have the same set 
of neighbors -namely, the central node- and
are therefore structurally equivalent (more precisely, they will be
false twins); then, after identifying all leaves in light of this 
twin notion for vertices,
the quotient graph amounts to $K_2$. 

But in the network literature
there is no equivalence notion for ``similar'' higher order subgraphs,
which would pave the way to a systematic reduction of 
(eventually defined) CSP structures.
As explained in detail in Section 
\ref{sec-bac} (see, specifically, subsection \ref{subsec-cp}), the goal of this paper is to fill this gap by introducing
a mathematical framework allowing
for a systematic classification of CSP networks and other 
partitioned structures. The key idea 
is to
introduce the concept of {\em twin subgraphs}, a notion
which extends to arbitrary
order that of twin (structurally equivalent) 
vertices. This mathematical framework will
be developed in Sections \ref{sec-falsetwins}
and \ref{sec-truetwins}, which address graph-theoretic problems
of independent interest (that is, problems which go beyond the eventual 
application of these notions to the classification of CSP structures).
These sections introduce and
elaborate on the idea of {\em F-twin} and {\em T-twin}
subgraphs, which in a sense are dual to each other and generalize 
several known properties of false twin and  true twin vertices; e.g.\ distinct connected components of
F-twin pairs will be proved to be disjoint and non-adjacent, whereas
disjoint T-twin pairs will be fully connected to each other.
With this background, the classification of CSP networks will then be 
tackled in 
Section \ref{sec-csp}. In Section \ref{sec-ex} we present the lines
along which these structures can be identified in real cases by
analyzing a subnetwork of the network of manufactures of metal arising
from 1994 world trade statistics. These data are available and 
analyzed in \cite{deNooy}, in the spirit of the
the aforementioned seminal work \cite{wallerstein74}, 
and nowadays
define a widely used benchmark for the positional analyses of networks.
Finally, Section \ref{sec-con} compiles
some 
lines for future research.


\section{Background on graphs, twins, and core-periphery networks}
\label{sec-bac}

\subsection{Graph-theoretic notions}
\label{subsec-graphs}

We refer the
reader to \cite{bollobas, bondy, diestel, harary69}
for excellent introductions
to graph theory. Throughout the paper we will work with
undirected graphs ${\cal G}=(V,E)$
without parallel edges or self-loops, so that edges 
can be thought of as pairs of distinct vertices (also termed
{\em nodes}). Given 
a graph ${\cal G}$, its vertex and edge sets will be 
written as $V({\cal G})$ and $E({\cal G})$, respectively, or
simply as $V$ and $E$ if there is no possible ambiguity.
We will only work with finite graphs, that is, the order
(number of vertices) will be finite in all cases. 
With 
notational 
abuse, we will often write $v \in {\cal G}$ to mean $v \in V({\cal G})$
and $V_0 \subseteq {\cal G}$ for $V_0 \subseteq V({\cal G})$.
Analogously, we will say that two graphs are disjoint when
their vertex sets are disjoint (note that the latter
implies that the edge sets are disjoint as well).

A {\em path} of length $k \geq 0$
is a graph with $k+1$ distinct vertices
$v_0,$ $v_1, \ldots, v_k$
and edges $e_1, \ldots, e_k$ with $e_i$ joining $v_{i-1}$ and $v_i$. 
Since we are not allowing parallel edges, a path is uniquely
defined by its vertex set. 
We say that $v_0$ and $v_k$ are {\em linked} by such a path.
When $k \geq 1$, sometimes the
vertex set
will be implicitly assumed to inherit the order defined by the
indices and we will then speak of a path {\em from} $v_0$ {\em to} 
$v_k$. 
The {\em distance}, $d$, between a pair of distinct vertices 
in the same connected component of a given graph
is the length of a shortest path linking them.
The {\em eccentricity} of a vertex in a connected graph is the maximum distance
to other vertices.
The distance between two disjoint 
subgraphs $H_1$ and $H_2$ 
lying in the same connected component of a given graph 
is defined as ${\rm min}  \{d(u,v), \ u \in H_1$, $v \in H_2\}$.
We say that two disjoint subgraphs $H_1$ and $H_2$
are not adjacent if there is no adjacent pair 
$(u,v)$  with $u \in H_1$, $v \in H_2$; if both subgraphs lie in the
same connected component of ${\cal G}$, this is equivalent
to saying that $d(H_1, H_2) \geq 2$.

We will denote by ${\cal N}(u)$ the set of neighbors of a 
given vertex $u$ (namely,
the set of vertices adjacent to $u$), 
and write ${\cal N}[u]={\cal N}(u) \cup \{u\}$.
The {\em degree} of a vertex $u$ is the number of elements
in ${\cal N}(u)$. We will call a vertex of degree one a {\em leaf}
(note that this term is often reserved to cases in which
the whole graph is acyclic, that is, a disjoint union of trees),
and will say that it is {\em attached} to its unique adjacent vertex.

The {\em null graph} 
defined by $V=\emptyset$ will be denoted by $K_0$; $K_n$ with $n \geq 1$
stands for the complete graph on $n$ vertices. 
The complement of a graph ${\cal G}=(V,E)$ of order $n$ 
(namely, $(V, E(K_n)-E)$) will be written as $\overline{\cal G}$, and
$E_n$ will stand for the {\em empty graph} $\overline{K_n}$ 
on $n \geq 1$
vertices.
Cycles, paths and stars on $n$ vertices 
will be written as $C_n$, 
$P_n$ and $S_n$, respectively, with $n \geq 3$ for cycles. As usual,
the union and intersection of 
${\cal G}_i=(V_i, E_i)$ 
($i=1,$ $2$) are the graphs $(V_1 \cup V_2, E_1 \cup E_2)$
and $(V_1 \cap V_2, E_1 \cap E_2)$, respectively.
The {\em join} ${\cal G}_1 + {\cal G}_2$ of two graphs with disjoint vertex sets
$V({\cal G}_1)$, $V({\cal G}_2)$
is the graph obtained after enlarging ${\cal G}_1 \cup {\cal G}_2$ 
with all possible edges joining
the vertices of ${\cal G}_1$ to those of ${\cal G}_2$
(sometimes we express the latter by saying 
that ${\cal G}_1$ and ${\cal G}_2$ are fully connected
to one another).



A {\em partitioned} graph is simply a graph 
whose vertex set is split into (pairwise disjoint) 
classes. A {\em $k$-partitioned graph} is a partitioned graph
with $k$ non-empty partition classes.
Obviously, a partitioned graph defines an equivalence relation in
the set of vertices. The {\em quotient graph} (often called a
{\em supergraph})
of a 
partitioned graph is defined as a graph whose
vertex set is the quotient set (that is,
vertices in the quotient graph correspond to the partition classes
in the original graph), two distinct
vertices in the quotient being adjacent if and only if
the original graph has at least one edge which joins vertices
belonging to the corresponding pair of classes. 

An {\em isomorphism} of two graphs ${\cal G}_1$ and ${\cal G}_2$
is a bijection $\varphi:V_1 \to V_2$ (with $V_i=V({\cal G}_i)$) which
preserves adjacencies, that is, such that 
any given pair of vertices $u$, $v$ 
in ${\cal G}_1$ are adjacent if and only if
$\varphi(u)$ and $\varphi(v)$ are adjacent in ${\cal G}_2$.
An isomorphism 
of partitioned
graphs is a graph isomorphism which keeps the classes invariant.


\subsection{Twins}  
\label{subsec-twins}


Different analytical and computational issues arise 
in connection 
to the 
existence and the distribution
of isomorphic copies 
of certain subgraphs of a given graph:
see e.g.\ 
\cite{alonBollobas, 
chungErdos, erdosHajnal, harary69, 
leeSIAM} 
and 
references therein.
From a different perspective, some attention has been focused on vertices which share
the same connection pattern within a graph. Such vertices receive (at least) two different names in the
literature, namely, {\em twins} and {\em structurally equivalent vertices}, as detailed in the
sequel. Two (distinct) 
vertices $u$ and $v$ are {\em false twins} (resp.\ {\em true twins}) if 
${\cal N}(u)={\cal N}(v)$  (resp.\ ${\cal N}[u]=
{\cal N}[v]$) \cite{bandelt, burlet, hernando2016, korach2008}. 
The exclusion of
self-loops yields $u \notin {\cal N}(u)$ and
this implies that false twins 
are not adjacent. In
the dual case, 
true twins are necessarily adjacent to each other:
for these reasons, true and false twins are also 
called {\em adjacent} and {\em non-adjacent twins}
(see e.g. \cite{bandelt, hernando2007, korach2008}). 
True twins correspond to 1-twins in the terminology of \cite{charonTwinFree, honkala}.
By contrast, in the social network analysis literature twin vertices $u$ and $v$ are said to
be {\em (weakly) structurally equivalent}: 
this means that the transposition $t_{u,v}$ of
$u$ and $v$ yields an automorphism of the graph 
(cf.\ \cite{borgatti94, brandes}), 
a condition
which is easily seen equivalent to $u$ and $v$ being (false or true) twins in the sense indicated
above. 
%

The F-twin and T-twin notions that will be
introduced in Sections \ref{sec-falsetwins} and 
\ref{sec-truetwins} for arbitrary subgraphs somehow
combine the two ideas at the beginning of the paragraph above. 
Twin subgraphs will
be isomorphic copies of each other and, additionally, they will share the connection
pattern to the remainder of the graph; in other words, our approach
will define a structural equivalence notion for (isomorphic) subgraphs 
which extends 
the one already defined for single vertices. Consistently, twin subgraphs will retain,
{\em mutatis mutandis,}
certain properties already known for twin vertices, such as the aforementioned
adjacency properties (which will hold for disjoint twin subgraphs;
cf.\
Corollaries \ref{coro-nonadj}
and \ref{coro-disjointTtwins}), the duality between F-twins and T-twins
in the sense that a pair of twins of one type defines a pair of the
other on the complement  graph
(Theorem \ref{th-duality}),
or the fact that twins will have the same distance multisets to 
the vertex set of the graph 
(cf.\ Proposition \ref{propo-distances}).
In particular, twin subgraphs will
define {\em homometric} sets (Corollary \ref{coro-homometric};
cf.\ \cite{albertsonHomometric,axenovichHomometric,rosen82}). 
Both notions will induce a classification
in the family of isomorphic copies of each induced subgraph, extending the
way in which false and true
twin concepts classify the vertices of a graph. These, together
with other related results, will be extensively discussed 
in Sections \ref{sec-falsetwins} and 
\ref{sec-truetwins}.


\subsection{\Cp 
networks}
\label{subsec-cp}

Consider one of the ``idealized'' \cp (CP) networks mentioned in Section \ref{sec-intr}, 
 namely, the one defined by a 2-partitioned graph 
with the following two classes of vertices:\vspace{-2mm}
\begin{itemize} 
\item[(i)] {\em core vertices}, \label{corevertices}which are fully connected to each other and also to 
the vertices in the second class (defined below);\vspace{-1mm}
\item [(ii)] {\em periphery vertices}, \label{peripheryvertices}totally disconnected from each other (and fully
connected to the core, in light of the first requirement above).
\end{itemize}
As indicated in the Introduction, other core-periphery connection patterns are possible, although
the one above is often used as a starting point in different analytical and computational approaches
to this topic (see e.g.\ \cite{borgatti99, deNooy}). 
These core-periphery networks are simply 2-partitioned graphs of the
form $K_p + E_r$ (find notations in subsection \ref{subsec-graphs};
when using a 2-partitioned structure in $K_p + E_r$,
we assume throughout the document and without further mention that
the two partition classes are the vertex sets of $K_p$ and $E_r$).
Cases with a unique core vertex amount to 
the star $S_n = K_{1,n-1}=K_1 + E_{n-1}$.
In the simplest setting ($n=2$) we get a 2-partitioned $S_2=K_2=K_1+K_1$, with a single core and a single periphery vertex; 
note that $E_1=K_1$, and we prefer to use
the latter notation for the singleton graph.


Aiming at later developments let us note that, in 
a certain sense, $K_2$ is substantially different from all other
joins $K_p+E_r$. Actually, we may think of $K_2=K_1+K_1$ as the quotient
graph of any other join of the form $K_p + E_r$. But, in order
to 
extend these ideas to support the definition and classification of
 more complex
structures, we emphasize that the reduction above 
comprises 
more than a quotient reduction. Indeed, all core vertices (namely, those
of $K_p$) are true twins as defined in 
subsection \ref{subsec-twins} above
and, 
analogously, all periphery vertices (the ones in $E_r$) are
false twins.
In this context, $K_2$ arises not only as the reduction
of other joins, but also as the unique {\em twin-free} network
meeting the requirements (i) and (ii) above. From this point
of view we may think of $K_2$ as the unique \cp {\em structure}
(we use the latter term to make a distinction with the CP networks 
$K_p + E_r$ above, which are allowed to display twin vertices).
To avoid any misunderstanding, let us clarify that $K_2$ is twin-free
only as a 2-partitioned graph, that is, we cannot consider
both vertices as (true) twins because they belong to different
partition classes;
cf.\ the beginning of Section \ref{sec-csp}.

However, when scaling these ideas to define formally \csp (CSP) 
structures, and eventually other structures with more partition classes,
one finds the problem that there is no appropriate
analog of the twin notions mentioned above for subgraphs with more
than one vertex. Since the intuitive idea behind the concept of 
a core is that of a set of heavily connected vertices, the true-twin
notion for single vertices may well apply to reduce the number
of admissible core subgraphs in these higher order
structures; by contrast, 
in the literature one finds no way to reduce conveniently
the semiperiphery-periphery 
subgraph. 

To put it in the simplest possible
setting, compare the CP network $K_1 + E_2$
(Fig.\ 1(a)), 
which amounts to a 2-partitioned path $P_3$ with one class (the core,
painted black in the figure)
defined by the central node, with a 3-partitioned path $P_5$ in which
the three classes are defined by the central vertex (core), the two
vertices with eccentricity three (semiperiphery vertices, grey) and the two
leaves (periphery vertices, white) (Fig.\ 1(b)). 
We may think of the latter as a (sometimes called) 
{\em spider graph} with a central
vertex (the core) and two legs,
each one a $P_2=K_2$ attached to the core by a single articulation
(the semiperiphery vertices). 

\begin{figure}[ht] \label{fig-cp2spider}
\mbox{} \vspace{3mm} \\
\mbox{} \hspace{32mm} 
\parbox{0.5in}{
\epsfig{figure=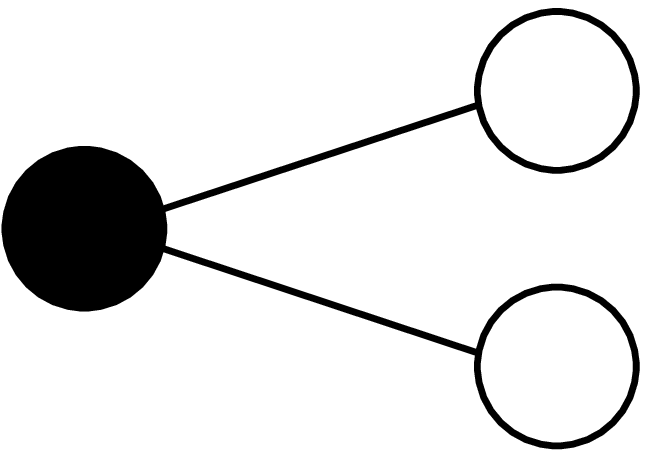, width=0.15\textwidth}
}
\hspace{37mm}
\parbox{0.5in}{
\epsfig{figure=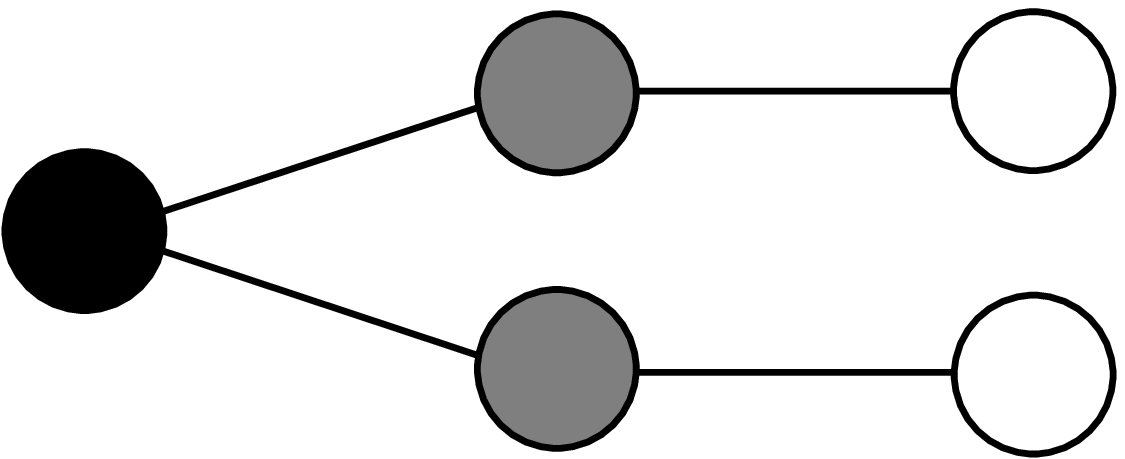, width=0.27\textwidth}
}
\vspace{4mm}
\caption{\ \ (a) CP network 
$K_1 + E_2$ 
\hspace{24mm} (b) A spider \hspace{35mm}}
\end{figure}

\vv

\noindent As indicated above,
in the CP case ($K_1 + E_2$) the false twin notion makes it
possible to identify the two peripheries into a single one,
reducing the network to a 2-partitioned $K_2=P_2$ (cf.\ Fig.\ 2(a)).
But, how can we reduce the CSP case (the spider) to a single $P_3$, which
captures the essential connection pattern? (Fig.\ 2(b)). Note
that both legs in Fig.\ 1(b) have exactly the same structure and, accordingly, 
we should find a systematic way to perform such reduction. Note also that
neither the semiperiphery vertices nor the periphery ones in Fig.\ 1(b)
are false twins, so that an eventual recourse to the notion
of twin vertices would fail for our present purpose.

\begin{figure}[ht] \label{fig-cpcsp}
\mbox{} \vspace{3mm} \\
\mbox{} \hspace{32mm} 
\parbox{0.5in}{
\epsfig{figure=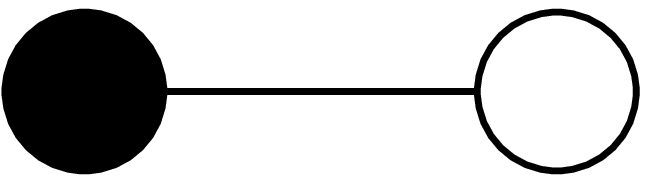, width=0.15\textwidth}
}
\hspace{37mm}
\parbox{0.5in}{
\epsfig{figure=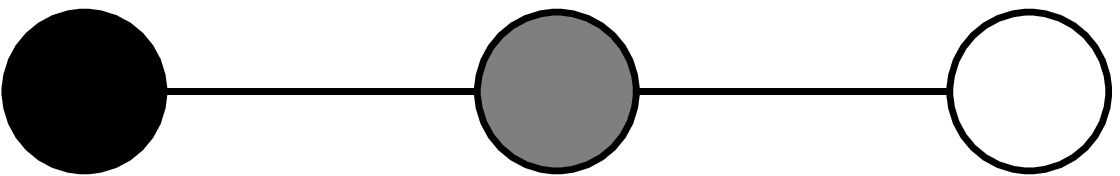, width=0.27\textwidth}
}
\vspace{4mm}
\caption{\ \ (a) CP structure
\hspace{30mm} (b) CSP structure \hspace{28mm}}
\end{figure}

\vv

\noindent
Obviously, it would be easy to identify equal-length legs in spider graphs;
however, more complex structures are possible: think e.g.\ of cases with more
cores and/or with other connection patterns within the semiperiphery
(actually, different CSP structures will arise in Sections \ref{sec-csp}
and \ref{sec-ex}; see Figs.\ 3-8).
Additionally,
the goal should be the development of a broader mathematical framework
allowing for an identification of (say) structurally equivalent,
higher order
subgraphs in greater generality. The idea is to formalize the notion
of isomorphic subgraphs or arbitrary order 
displaying, in a sense to be made precise, 
the same connection patterns to the remainder
of the graph, generalizing the false-twin and true-twin concepts for
single vertices. 
The F-twin notion for arbitrary subgraphs, 
together with the dual concept
of T-twin subgraphs, are 
aimed at filling this gap. After introducing and 
discussing these ideas in Sections \ref{sec-falsetwins}
and \ref{sec-truetwins}, we will be back to 
CSP structures 
in Sections \ref{sec-csp} and \ref{sec-ex}.


\section{F-twin subgraphs}
\label{sec-falsetwins}





\subsection{Definition and elementary properties}

\begin{defin} \label{defin-Ftwins}
Let $H_1$ and $H_2$ be two 
induced subgraphs of a 
graph ${\cal G}$. Denote by
$V_i$ the vertex set of $H_i$. 
$H_1$ and $H_2$ are called 
{\em F-twins} if
they are isomorphic via a map $\varphi:V_1 \to V_2$ for which the identities
\begin{equation} \label{Ns}
{\cal N}(u)- V_1 = {\cal N}(\varphi(u))- V_2
\end{equation}
hold for all $u \in V_1.$
\end{defin}

\vv

\noindent We may also say that the set of vertices $V_1$ and $V_2$ are F-twins, since the
definition above requires $H_1$ and $H_2$ to be the subgraphs
induced by $V_1$ and $V_2$ and there is no possible ambiguity.
The reason for the requirement that F-twins 
are induced subgraphs 
should become apparent
in light of a simple example, defined by 
the graph ${\cal G}=P_3 \cup K_3$. Let $H_1$ be the $P_3$-component 
of ${\cal G}$
and $H_2$ any one of the three subgraphs of $K_3$ isomorphic to $P_3$.
Should F-twins not be required to be induced subgraphs,
$H_1$ and $H_2$ would be F-twins, because the identities
(\ref{Ns}) hold trivially
since both sides are empty for all vertices. However there exists
an extra edge in $K_3$ which make the endvertices of $H_2$ adjacent 
in ${\cal G}$ without the endvertices of $H_1$ being so. 
Since the idea of the F-twin notion is to capture
identical adjacency patterns, we rule out this type
of situations by requiring
$H_1$ and $H_2$ to be induced subgraphs.

Note that any induced subgraph is trivially an F-twin of itself; 
we will say that a given induced subgraph is a {\em proper F-twin}
if it has at least an F-twin different 
from itself (and both of them will also be said to be proper F-twins
of each other). A {\em trivial} F-twin is an induced subgraph that has no
F-twin but itself.

In particular, the notion above for two single distinct vertices
$u$, $v$ amounts to requiring that they are false twins
in the sense that
${\cal N}(u)={\cal N}(v)$, as defined in subsection \ref{subsec-twins}.
Just note that 
$u \not\in {\cal N}(u)$ and $v \not\in {\cal N}(v)$, so that
(\ref{Ns}) holds in this case if and only if ${\cal N}(u)={\cal N}(v)$.

\vv

\begin{propo} \label{propo-nonconnected}
Two 
induced subgraphs $H_1$ and $H_2$ are F-twins if and only if
their 
connected components can be matched 
as pairs of F-twins.
\end{propo}

\noindent {\bf Proof.} Assume first that $H_1$ and $H_2$ are 
F-twins, and let $\varphi$ denote the isomorphism arising in Definition
\ref{defin-Ftwins}. Then $\varphi$ induces $k$ isomorphisms $\varphi_1,\
\ldots, \varphi_k$ between the connected components of $H_1$ and $H_2$;
denote these connected components by 
$H_{i,j} $, with $i \in \{1,2\}, \ j \in \{1, \ldots, k\}$, and let
accordingly $V_{i,j} $ be the
vertex set of  $H_{i,j}$, so that $\varphi_j:V_{1,j}\to V_{2,j}$. Then obviously
$$V_i = \bigcup_{j=1}^k V_{i,j}$$ 
and, provided that a vertex $u$ (resp.\ $\varphi(u)$) 
belongs to $V_{1,j}$ (resp.\ to $V_{2,j}$), it is also clear that
${\cal N}(u) \cap V_1 \subseteq V_{1,j}$, 
which implies
${\cal N}(u) \cap V_{1,k}= \emptyset$ if $k \neq j$
(resp.\ ${\cal N}(\varphi(u)) \cap V_2
\subseteq V_{2,j}$ and then
${\cal N}(\varphi(u)) \cap V_{2,k}= \emptyset$ if $k \neq j$). This yields
\begin{equation} \label{Nsconn}
{\cal N}(u)- V_{1,j} = {\cal N}(u)- V_1 = 
{\cal N}(\varphi(u))- V_2 = {\cal N}(\varphi_j(u))- V_{2,j}
\end{equation}
so that $H_{1,j}$ and $H_{2,j}$ are  indeed F-twins.

The converse result proceeds in exactly the same manner and details
are left to the reader.\vspace{-3mm} 

\hfill $\Box$

\noindent The result above is non-trivial only when $H_1$ and $H_2$ are
not connected. In this setting, even if $H_1$ and $H_2$ are proper F-twins some
of their components might be trivial F-twins. 

\vv

\begin{propo} \label{propo-disjoint} 
If $H_1$ and $H_2$ are 
proper F-twins, 
the intersection
$V_1 \cap V_2$, if non-empty, induces 
a set 
of connected components of both $H_1$ and $H_2$. 
\end{propo}

\noindent {\bf Proof.} Assume that $w \in V_1 \cap V_2$, 
and let 
$K_1$ and $K_2$ be the connected components 
of $H_1$ and $H_2$ which accommodate $w$. 
Assume that
 $K_1 \neq K_2$ 
and, w.l.o.g., suppose that there is a vertex in $K_2$ not
belonging to $K_1$.
The set $V(K_2)$ can be described as the disjoint union of
$V(K_1) \cap V(K_2)$ and 
$V(K_2) - V(K_1)$ and,
since $K_2$ is connected,
there must exist two adjacent vertices $u$, 
$v$ with $u \in V(K_1) \cap V(K_2)$ and $v \in V(K_2) - V(K_1)$.
The fact that $v \notin V(K_1)$ implies $v \notin V_1=V(H_1)$; indeed,
should it belong to $V_1$, since it is adjacent to $u \in V(K_1) \subseteq V_1$ 
it would 
necessarily belong to the same connected component of $u$, that
is, to $K_1$, but we know that $v \notin V(K_1)$.
This implies that $v \in {\cal N}(u)-V_1$ and, in light of (\ref{Ns}), it must happen that
$$v \in {\cal N}(\varphi(u)) - V_2,$$
whoever $\varphi(u)$ is.
But this is impossible because $v \in V(K_2) \subseteq V_2$ implies $v \notin {\cal N}(\varphi(u)) - V_2$.
Hence $K_1=K_2$ and since $V(K_1) \subseteq V_1$, $V(K_2) \subseteq V_2$
we conclude that the whole connected components $K_1=K_2$ are in 
the intersection $H_1 \cap H_2$ as we aimed to show.\vspace{-3mm}

\hfill $\Box$

\noindent In particular, Proposition \ref{propo-disjoint}
implies that distinct connected F-twins
are actually disjoint. 

\vs

\begin{coro} \label{coro-disjoint} 
If $H_1$ and $H_2$ are connected proper F-twins
then $V_1 \cap V_2 = \emptyset$. 
\end{coro}







\subsection{Distance-related properties}

We know from Proposition \ref{propo-disjoint}
that non-empty intersections of F-twins necessarily 
span connected components of both. 
On the other hand, when 
two F-twin subgraphs are disjoint one can easily show 
that they cannot be adjacent (just derive from (\ref{Ns})
the identities ${\cal N}(u) \cap V_2 = \emptyset$
for all $u \in V_1$).
A stronger statement actually holds. 

\begin{propo} \label{propo-distnueva}
If $H_1$ and $H_2$ are disjoint F-twins in a given graph ${\cal G}$,
 then 
their connected components  
can be arranged as F-twin pairs $(H_{1,j}, H_{2,j})$
in a way such that, for every $j$, \vspace{-3mm}
\begin{itemize}
\item either $H_{1,j}$ and $H_{2,j}$
are connected components of ${\cal G}$;  or \vspace{-1mm} 
\item both $H_{1,j}$ and $H_{2,j}$
belong to the same connected component of ${\cal G}$ and
$d (H_{1,j}, H_{2,j})=2$. 
\end{itemize}
\end{propo}





\noindent {\bf Proof.} 
Take a connected component $H_{1,j}$ of $H_1$  and
assume that there 
exists a vertex $v \notin V_1=V(H_1)$ adjacent to some $u \in V_{1,j}=
V(H_{1,j})$. 
In light of (\ref{Ns}), it follows that 
$v \in {\cal N} (\varphi(u))-V_2$, with $V_2=V(H_2)$; this implies 
that $v \notin V_2$ (a property that will be used later)
and also that $(u, v, \varphi(u))$ is a path. Let $H_{2,j}$ 
be the connected component of $H_2$ accommodating $\varphi(u)$: then
$H_{1,j}$ and $H_{2,j}$  are isomorphic via $\varphi$; moreover,
they are in the same connected component of ${\cal G}$ and,
additionally, $d(H_{1,j}, H_{2,j}) \leq 2$. The aforementioned property
that any vertex $v \notin V_1$ adjacent to $u \in H_{1,j}$ cannot
belong to $V_2=V(H_2)$ shows that, actually, $d(H_{1,j}, H_{2,j}) = 2$.


The same reasoning applies to all connected components of $H_1$.
Those for which there is no adjacent vertex away from
$V(H_1)$ are by definition connected components of ${\cal G}$.
Exactly the same reasoning applies to the connected components
of $H_2$ and this completes the proof.


\hfill $\Box$


\noindent Note also that 
for components $H_{1,j}$, $H_{2,k}$ of $H_1$ and $H_2$ which
do not define an F-twin pair and which are
contained in the same connected component 
of ${\cal G}$  it holds as well
that $d(H_{1,j},H_{2,k}) \geq 2$
since they cannot be adjacent to each other.


Corollary \ref{coro-nonadj} follows directly
from Proposition \ref{propo-distnueva}. 
Implicit in its first claim
is the fact that connected, proper F-twins which are not connected component
themselves must lie in the same connected component of ${\cal G}$.
The second claim emphasizes 
that our notion extends the non-adjacency property 
of false twin vertices mentioned in subsection \ref{subsec-twins}.

\vv

\begin{coro} \label{coro-nonadj} 
If $H_1$ and $H_2$ are connected proper F-twins in a given graph ${\cal G}$, 
then either they are 
connected components of ${\cal G}$ or 
$d (H_1, H_2)=2$. 
In either case, connected proper F-twins are not adjacent to each other.
\end{coro}

 
Another distance-related property of proper F-twins is that
they are {\em homometric}; this means that
the  distance multisets of both are the same \cite{albertsonHomometric,
axenovichHomometric, rosen82}. 
The distance multiset of
an order-$k$ subgraph $H$ of a connected graph ${\cal G}$ is the multiset of
$\binom{k}{2}$
distances (in 
${\cal G}$) between vertices of $H$.

\begin{lema} \label{lema-mirrorpaths}
Assume that $H_1$ and $H_2$ are disjoint 
F-twin subgraphs
of a 
graph ${\cal G}$. Let $(u_0, \ldots, u_k)$ be a vertex
sequence defining a path (of length $k$)
in ${\cal G}$. Then $(v_0, \ldots, v_k)$, with
$$v_i = \begin{cases} \varphi(u_i) & \text{ if } u_i \in V(H_1)\\
  \varphi^{-1}(u_i) & \text{ if } u_i \in V(H_2) \\
u_i & \text{ if } u_i \notin V(H_1) \cup V(H_2), \\ 
\end{cases}$$
also defines a length-$k$ path.
\end{lema}

\noindent {\bf Proof.} The fact that all vertices $v_i$ are distinct
is a direct consequence of the construction: indeed, note that
$\varphi$ maps $V_1=V(H_1)$ onto $V_2=V(H_2)$ and, conversely,
$\varphi^{-1}$ maps $V_2$ onto $V_1$. Since $V_1$, $V_2$ and
$V-(V_1 \cup V_2)$ (with $V=V({\cal G})$) 
are pairwise disjoint sets, then the claim follows easily from the facts
that $\varphi$, $\varphi^{-1}$ and the identity are bijections
and that the vertices $u_i$ are all distinct.

The other fact that needs to be proved is that the pairs $\{v_{i-1}, v_i\}$
are adjacent. Since we know that disjoint F-twins are not adjacent
(cf.\ Proposition \ref{propo-distnueva} and Corollary \ref{coro-nonadj})
and the isomorphisms $\varphi$ and $\varphi^{-1}$ preserve adjacencies,
we only need to check that $v_{i-1}$ and $v_i$ are adjacent when
one of them (say $v_{i-1}$) belongs to one of the twins (e.g.\ to $H_2$,
for later notational simplicity)
and $v_i$ is not in $H_1 \cup H_2$. This means that $v_{i-1} = \varphi(u_{i-1})$
with $u_{i-1} \in V_1$ and that $v_i = u_i \notin V_1 \cup V_2.$
Now use the fact that $u_i \in {\cal N}(u_{i-1})$ because the vertices
$u_i$ define a path. Additionally, since $u_i \notin V_1$, from (\ref{Ns})
we conclude that $u_i \in {\cal N}(\varphi(u_{i-1})) - V_2$.
The identities $v_i = u_i$, $v_{i-1}=\varphi(u_{i-1})$, 
show that $v_i \in {\cal N}(v_{i-1})$, as we aimed to prove. 

\hfill $\Box$

\begin{propo} \label{propo-distances}
Assume that $H_1$ and $H_2$ are disjoint  F-twin subgraphs
of a connected graph ${\cal G}$, and let $u \in V(H_1)$. Then, for any 
other vertex $\tilde{u}$ in 
${\cal G}$ the following assertions hold.
\vspace{-3mm}
\begin{itemize}
\item[a)] If $\tilde{u} \in V(H_1)$, then $d(u,\tilde{u})=d(\varphi(u), \varphi(\tilde{u}))$.
\vspace{-1mm} 
\item[b)] If $\tilde{u} \in V(H_2)$, then $d(u,\tilde{u})=d(\varphi(u), \varphi^{-1}(\tilde{u}))$.
\vspace{-1mm}
\item[c)] If $\tilde{u} \notin V(H_1) \cup V(H_2)$, then 
$d(u,\tilde{u})=d(\varphi(u), \tilde{u})$.
\vspace{-1mm}
\end{itemize}
\end{propo}

\noindent {\bf Proof.} The results follow in a straightforward
manner from Lemma \ref{lema-mirrorpaths} since the
set of paths from $u$ to $\tilde{u}$ are in a one-to-one, length-preserving 
correspondence to the ones that link $\varphi(u)$
to $\varphi(\tilde{u})$, $\varphi^{-1}(\tilde{u})$ or $\tilde{u}$,
depending on the case. The distance identities follow as an
immediate consequence simply because the distance between two vertices is 
the minimum length of the paths linking those vertices.

\hfill $\Box$

\noindent Another way to state item a) of Proposition \ref{propo-distances}
is the following.

\vv

\begin{coro} \label{coro-homometric}
Disjoint F-twin subgraphs
of a connected graph ${\cal G}$ are homometric.
\end{coro}

\noindent Note also that c) extends a known property of false twin vertices
(cf. \cite[Proposition 1.1]
{hernando2016}).


\subsection{On the classification of F-twin subgraphs}
\label{subsec-classification}



The F-twin notion classifies the set of isomorphic copies of 
any induced subgraph of a given graph, as shown below.



\begin{theor} \label{th-equiv}
Let $H$ be an induced subgraph of ${\cal G}$ and denote by
${\cal H}$ the set of induced subgraphs of ${\cal G}$ which
are isomorphic to $H$. Then the F-twin relation stated in 
Definition \ref{defin-Ftwins}
is an equivalence relation in ${\cal H}$.
\end{theor}

\noindent {\bf Proof.} The F-twin relation is obviously 
reflexive since we may
set $\varphi$ as the identity when $H_1=H_2$
in Definition \ref{defin-Ftwins}.
The fact that it is also symmetric is also easily checked, just
using the inverse $\varphi^{-1}$ of the
isomorphism $\varphi$.  
Transitivity is also rather straightforward.
Let us assume that ($H_1$, $H_2$)
and ($H_2$, $H_3$) are pairs of F-twins,
and denote by $\varphi$ and $\psi$ the
isomorphisms between $H_1$ and $H_2$ and between $H_2$ and $H_3$,
respectively. One can check that 
the isomorphism $\zeta = \psi \circ \varphi$ yields
\begin{equation} \label{Nsbis}
{\cal N}(u)- V_1 = {\cal N}(\zeta(u))- V_3
\end{equation}
for all $u \in V_1$: indeed, this is an immediate 
consequence of (\ref{Ns}) 
and the corresponding identity for the isomorphism $\psi$, that is,
${\cal N}(v)- V_2 = {\cal N}(\psi(v))- V_3$
for all $v \in V_2$. The identities (\ref{Nsbis}) are obtained 
just by setting 
$v=\varphi(u)$.

\hfill $\Box$

\noindent Since all these classifications of induced subgraphs 
eventually act on the same underlying object (the graph itself), it is
natural to wonder about possible interrelations between
such classifications of different subgraph families. In the forthcoming
subsections we provide some initial results in this direction;
we explore, in particular, whether F-twin vertices may belong to larger
connected F-twin structures, and also provide some remarks about the
F-twin classification of the family (to be denoted as ${\cal H}_2$) 
of subgraphs isomorphic
to $K_2$. 
With terminological
abuse we will refer to this problem as the classification of F-twin
edges (namely, we deliberately identify 
an edge $e$ with the $K_2$-graph induced by its endvertices
$u, v$, the latter being in fact the graph $\left(\{u,v\}, \{e\}\right)$):
with this cautionary remark in mind
the reader can think of ${\cal H}_2$ simply as the set of edges.


\subsubsection{F-twin vertices within larger F-twin structures}

Assume that a given graph has 
a class of three or more F-twin vertices. We know that they
are pairwise non-adjacent and, by definition, that they share
a common set of neighbors. It then follows that any two proper subsets
 of this class with the same number of elements (which
induce two empty graphs with the same number of vertices) are 
themselves F-twins, since any isomorphism matching 
the vertices of these two empty graphs preserves the 
relations involved in (\ref{Ns}). 
The other way round, we may think of this as an example
in which two proper F-twin subgraphs contain 
two proper F-twin vertices (more precisely, in a way such
that each vertex lies on one of the larger
twins), consistently with Proposition \ref{propo-nonconnected}. 
As shown below, this cannot happen, however, if such an F-twin vertex
is adjacent to at least another vertex in the larger twin; this
essentially means that the inclusion of pairs of F-twin vertices
into pairs of larger F-twin structures
is specific to singletons of these larger subgraphs.




\vv

\begin{propo} \label{propo-sub}
Assume that $u$ and $\varphi(u)$ are proper F-twin vertices. If $u$ is
properly 
contained in a connected 
proper F-twin $H$, 
then the F-twin vertex $\varphi(u)$ also belongs to $H$.
\end{propo}

\noindent {\bf Proof.} Let $v$ be a vertex in $H$ adjacent to $u$; such a vertex is guaranteed
to exist because $u$ is assumed to be 
properly contained 
in the connected 
subgraph $H$. The F-twin vertices $u$ and $\varphi(u)$ 
are known to verify the relation
${\cal N}(u) = {\cal N}(\varphi(u))$,
and $v \in {\cal N}(u)$ 
then yields $v \in {\cal N}(\varphi(u))$; for later use we recast this relation as
$\varphi(u) \in {\cal N}(v)$.


Let us suppose that $\varphi(u) \notin V(H)$, and denote by $\psi$ the isomorphism mapping
$H$ to its F-twin $\psi(H)$. For this F-twin relation, the identities (\ref{Ns}) yield
in particular for $v \in H$ 
$${\cal N}(v)-V(H) = {\cal N}(\psi(v))-\psi(V(H)).$$ 
Now, if $\varphi(u) \notin V(H)$ and given the fact that $\varphi(u) \in {\cal N}(v)$ as shown
above, we obtain $\varphi(u) \in {\cal N}(\psi(v))$; as before, we recast this as
$\psi(v) \in {\cal N}(\varphi(u))$. But using again 
${\cal N}(u) = {\cal N}(\varphi(u))$ we would get
$\psi(v) \in {\cal N}(u)$ and this is in contradiction with Corollary \ref{coro-nonadj}
because $u \in H$ and $\psi(v) \in \psi(H)$, meaning that the 
connected F-twin structures
$H$ and $\psi(H)$ would be adjacent to each other. This implies that necessarily
 $\varphi(u) \in V(H)$ and the claim is proved.

\hfill $\Box$



\noindent Corollary \ref{coroK1K2} follows from the case in which the proper
F-twin $H$ in Proposition \ref{propo-sub}
is isomorphic to $K_2$. In this case there is no way in which $H$ may
accommodate two distinct F-twin vertices, since they would obviously be 
adjacent to each other and this would contradict 
Corollary \ref{coro-nonadj}.

\vv

\begin{coro} \label{coroK1K2}
Vertices and edges admitting proper F-twins define mutually disjoint vertex sets.
\end{coro}







\noindent We finish this section with a pretty obvious but useful remark following Corollary \ref{coroK1K2}.

\vv

\begin{coro}
Graphs of order $\leq 5$ cannot display simultaneously proper F-twin vertices and  proper F-twin edges.
\end{coro}


\subsubsection{Non-trivial vertex set intersections between classes of F-twin edges}


Obviously, in any graph the classification of F-twin vertices
yields pairwise disjoint vertex classes.
Things may get more involved when studying the interrelation between 
different F-twin classes of subgraphs not isomorphic to a single
vertex. For instance, a 6-cycle (cf.\ the proof of Proposition
\ref{propo-c6} below) accommodates three pairs of F-twin edges
with non-empty vertex intersections among classes. In a way,
such a 6-cycle is the essential structure to signal this 
phenomenon. We recall that ${\cal H}_2$ denotes the set of subgraphs of ${\cal G}$ isomorphic 
to $K_2$.

\vv

\begin{propo} \label{propo-c6}
Assume that two elements of ${\cal H}_2$ within a graph ${\cal G}$
belong to different proper F-twin classes and have a 
common vertex. Then ${\cal G}$ contains
the cycle $C_6$ as an induced subgraph.
\end{propo}

\noindent {\bf Proof.} Let $H_1$ and $J_1$ be two subgraphs 
in ${\cal H}_2$ (namely, isomorphic to $K_2$) which belong to different nontrivial F-twin classes,
and denote by $H_2$ and $J_2$ two proper F-twins of $H_1$ and $J_1$, respectively
(with the corresponding isomorphisms to be denoted by $\varphi$ and $\psi$).
Assume that $v$ belongs to both $H_1$ and $J_1$, and let $u$ and $w$
be the other vertex of $H_1$ and $J_1$, respectively. 
We claim that $\varphi(u)=\psi(w)$ and that
the subgraph induced by $\{u, v, w, \varphi(v), \varphi(u), \psi(v)\}$
is a 6-cycle.

To show this, write the F-twin identity for $v \in H_1$ 
as
\begin{equation} \label{vH}
{\cal N}(v)-V(H_1) = {\cal N}(\varphi(v))-V(H_2).
\end{equation}
Since $w \in {\cal N}(v)$ and $w \notin V(H_1)$, we derive 
\begin{equation} \label{aux4}
w \in {\cal N}(\varphi(v)).
\end{equation}
For later use, notice that this implies $\varphi(v) \notin V(J_2)$ 
(that is, $\psi(v) \neq \varphi(v)  \neq \psi(w)$), 
since otherwise there would
be two adjacent vertices in $J_1$ and $J_2$ (namely, $w$ and $\varphi(v)$),
against Corollary  \ref{coro-nonadj}.

Note that $v$ also belongs to $J_1$ and therefore, analogously, 
${\cal N}(v)-V(J_1) = {\cal N}(\psi(v))-V(J_2)$
and, proceeding as above (use $u \in {\cal N}(v)- V(J_1)$), we get
\begin{equation} \label{aux4bis}
u \in {\cal N}(\psi(v)),
\end{equation}
and also $\psi(v) \notin V(H_2)$, that is $\psi(v) \neq \varphi(u)$ (we already
knew that $\psi(v) \neq \varphi(v)$). 

Now, restate (\ref{aux4}) as $\varphi(v) \in {\cal N}(w)$ and, from the fact that
$\varphi(v)\notin V(J_1)$ (to check this just 
note that $v \neq \varphi(v) \neq w$, 
the latter being clear in the light of (\ref{aux4})) and
the F-twin identity for $w \in V(J_1)$,
\begin{equation} \label{wJ}
{\cal N}(w)-V(J_1) = {\cal N}(\psi(w))-V(J_2),
\end{equation}
derive $\varphi(v) \in {\cal N}(\psi(w))$
or, equivalently,
\begin{equation} \label{aux5}
 \psi(w) \in {\cal N}(\varphi(v)).
\end{equation}

We show in the sequel that, indeed, it is $\psi(w)=\varphi(u)$. Suppose $\psi(w)\neq \varphi(u)$;
as shown above we have $\psi(w) \neq \varphi (v)$ 
and both conditions
together would mean $\psi(w)  \notin V(H_2)$. 
Equations (\ref{vH}) and (\ref{aux5})
would then yield $\psi(w) \in {\cal N}(v)$. But then $v \in V(J_1)$ and 
$\psi(w) \in V(J_2)$ would
be adjacent to each other. We conclude that necessarily $\psi(w)=\varphi(u)$, as claimed.

The fact that $u, v, w, \varphi(v), \varphi(u)=\psi(w), \psi(v)$ yield a 6-cycle
follows from the adjacency relations defined by $H_1$, $J_1$, (\ref{aux4}), $H_2$, $J_2$, and (\ref{aux4bis}),
respectively.
It only remains to show that this cycle is actually induced by these vertices, namely, that
there are no additional adjacencies among them. Apart from the six edges defining the aforementioned
cycle, there are other nine possible links between the six vertices listed above; seven of these
are ruled out by Corollary \ref{coro-nonadj} (namely, those connecting $u,$ $v$ with $\varphi(u),$
$\varphi(v)$, since both pairs define the F-twins $H_1$, $H_2$, respectively,
and $v,$ $w$ with $\psi(v)$, $\psi(w)$,
which define $J_1$ and $J_2$; note that $\varphi(u)=\psi(w)$ and therefore
the pairs $\{v, \varphi(u)\}$ and $\{v, \psi(w)\}$ are the same). The two remaining pairs
are $\{u, w\}$ and $\{\varphi(v),\psi(v)\}$; consider the first one and note that $u \notin V(J_1)$,
so that the assumption $u \in {\cal N}(w)$ would imply $u \in {\cal N}(\psi(w))$
in light of (\ref{wJ}), but this is impossible because $u \in V(H_1)$ and $\psi(w)=\varphi(u) \in V(H_2)$
cannot be adjacent to each other. The fact that $\varphi(v)$ cannot be adjacent to $\psi(v)$ 
can be checked in the same terms, 
and the proof is complete.

\hfill $\Box$

\noindent We close this section by saying that the classification of
F-twin structures (beyond F-twin vertices) possibly defines 
other mathematical problems of interest. This is a topic for future study.


\section{T-twins}
\label{sec-truetwins}

We present in this section the dual concept of 
T-twin subgraphs, which extends the notion of true twin vertices 
discussed in subsection \ref{subsec-twins}.
This section will be briefer than the previous one; we just aim
at providing a complete framework extending to arbitrary subgraphs
the idea behind false and true twin vertices. We will also show 
(Theorem \ref{th-duality}) that in a precise sense 
the notions supporting F-twins and T-twins
are dual to each other, again extending a known property of
false and true twin vertices \cite{brandes, hernando2016}.

\vs

\begin{defin} \label{defin-Ttwins}
Let $H_1$ and $H_2$ be two 
induced subgraphs of a graph ${\cal G}$ and denote by
$V_i$ the vertex set of $H_i$. 
$H_1$ and $H_2$ are called 
{\em T-twins} if
they are isomorphic via a map $\varphi:V_1 \to V_2$ for which the identities
\begin{equation} \label{Ns-T}
{\cal N}(u) \cup V_1 = {\cal N}(\varphi(u)) \cup V_2
\end{equation}
hold for all $u \in V_1.$
\end{defin}

\noindent Again this extends the notion of true twin vertices introduced in subsection
\ref{subsec-twins}, which are defined by the identities
${\cal N}[u]={\cal N}[v]$, that is,
${\cal N}(u) \cup \{u\}={\cal N}(v) \cup \{v\}$, consistently with
(\ref{Ns-T}).

%
%
%
%
%

As in the F-twin case, 
we use the term {\em proper} T-twins for distinct T-twins.

\begin{propo} \label{propo-fully}
Let $H_1$ and $H_2$ be T-twins.
Then 
$V_1 \cap V_2$, 
$V_1-V_2$ and
$V_2 - V_1$ are fully connected to each other. 
\end{propo}

\noindent {\bf Proof.} From (\ref{Ns-T}) it is clear that
all vertices in $V_2 - V_1$ belong to ${\cal N}(u)$ for all $u \in V_1$,
and this means that $V_2 - V_1$ is fully connected to $V_1$ (in particular,
to $V_1-V_2$). Analogously,
$V_1 - V_2$ is fully connected to $V_2$. Using both properties together
we conclude that the intersection $V_1 \cap V_2$ is fully connected
to both $V_1-V_2$ and $V_2-V_1$ and the claim is proved. 

\hfill $\Box$

\begin{coro} \label{coro-disjointTtwins}
If $H_1$ and $H_2$ are disjoint T-twins, then $V_1$ is fully connected to 
$V_2$.
\end{coro}

\noindent The following result gives a precise meaning to the claim that
the F-twin and T-twin notions are dual to each other.

\vs

\begin{theor} \label{th-duality}
Two induced subgraphs $H_1$ and $H_2$ of a given graph ${\cal G}$
are T-twins (resp.\ F-twins) if and only if $\overline{H_1}$ and
$\overline{H_2}$ are F-twins (resp.\ T-twins) in $\overline{{\cal G}}$.
\end{theor}



\noindent {\bf Proof.} 
The reader can check in advance that if $H$ is an induced subgraph
of ${\cal G}$, then $\overline{H}$ is an induced subgraph
of $\overline{\cal G}$.
Assume now that $H_1$ and $H_2$ are T-twins, and let $\varphi$
be the isomorphism arising in Definition \ref{defin-Ttwins}; 
one can see that $\varphi$ is also an isomorphism between
the complements $\overline{H_1}$ and $\overline{H_2}$.
Denoting by $\overline{{\cal N}}(u)$ the neighborhood of $u$
in $\overline{{\cal G}}$, we need to show that the identities 
\begin{equation} \label{goal}
\overline{{\cal N}}(u)- V_1 = 
\overline{{\cal N}}(\varphi(u))- V_2 
\end{equation}
hold in $\overline{{\cal G}}$ for all 
$u$ in $V_1=V(\overline{H_1})=V(H_1)$.
We use the fact that 
$$\overline{{\cal N}}(u)=V({\cal G})-({\cal N}(u) \cup \{u\}),
\
\overline{{\cal N}}(\varphi(u))=V({\cal G})-({\cal N}(\varphi(u)) \cup \{\varphi(u)\})
$$
by definition of the complement.
These relations yield
\begin{equation} \label{aux2}
\overline{{\cal N}}(u)-V_1=V({\cal G})-({\cal N}(u) \cup \{u\} \cup V_1)=
V({\cal G})-({\cal N}(u) \cup V_1)
\end{equation}
(where we have used $u \in V_1$)
and, analogously,
\begin{equation} \label{aux3}
\overline{{\cal N}}(\varphi(u))-V_2=
V({\cal G})-({\cal N}(\varphi(u)) \cup \{\varphi(u)\} \cup V_2)=
V({\cal G})-({\cal N}(\varphi(u)) \cup V_2).
\end{equation}
The relations depicted in (\ref{goal}) then follow from (\ref{aux2})
and (\ref{aux3}) because 
$H_1$ and $H_2$ are T-twins, which means ${\cal N}(u) \cup V_1 = {\cal N}(\varphi(u)) \cup V_2$.

Both the case in which $H_1$ and $H_2$ are F-twins and the converse
results proceed in the same manner and details are left to the reader.

\hfill $\Box$

\noindent At first sight, a reader might
be slightly surprised with Theorem \ref{th-duality} since
T-twins may have non-empty intersections in the vertex sets and
(connected proper) F-twins seemingly not, as stated in Corollary
\ref{coro-disjoint}. But note that the latter holds 
as a consequence of Proposition \ref{propo-disjoint}
for {\em connected} F-twins: now assume $V_1 \cap V_2 \neq \emptyset$
for (even possibly connected) T-twins $H_1$, $H_2$.
From Proposition \ref{propo-fully} it follows that $V_1 \cap V_2$
is fully connected to both $V_1 - V_2$ and to $V_2 - V_1$, so that,
in the complementary (F-twin) subgraphs $\overline{H_1}$ and $\overline{H_2}$,
$V_1 \cap V_2$ is isolated from both $V_1 - V_2$ and $V_2 - V_1$.
This means that $V_1 \cap V_2$ induces a set of connected
components of both $\overline{H_1}$ and $\overline{H_2}$
and there is no contradiction with Proposition \ref{propo-disjoint}.

Finally, we mention that the T-twin relation also induces a classification
in the families ${\cal H}$ of isomorphic copies of induced subgraphs
$H$. Details are entirely analogous to those in 
Theorem \ref{th-equiv}  and are left to the reader.

\section{\Csp structures}
\label{sec-csp}

We take now a look back at 
subsection \ref{subsec-cp}; specifically, we provide here
a definition of \csp (CSP) structures extending the ideas presented there
and reducing the number of structures via the exclusion
of twin substructures,
according to the notions introduced
in Sections \ref{sec-falsetwins} and \ref{sec-truetwins}. We will work
in this section with 3-partitioned graphs (cf.\ subsection \ref{subsec-graphs})
and we make the remark that the F-twin and T-twin
notions introduced in Definitions \ref{defin-Ftwins} and \ref{defin-Ttwins}
apply also in this context just by assuming that
the isomorphism $\varphi$ is now an isomorphism of partitioned graphs,
namely, that it leaves the classes invariant (it maps core vertices into
core vertices, etc.).

\subsection{A parameterized definition of \csp structures}

We first note that the condition 
depicted in item (i) on page \pageref{corevertices},
defining core vertices, may be recast as the requirement that all of them
have eccentricity one.
This approach is intimately related to the {\em closeness centrality} notion,
widely used in network theory \cite{brandes, newman}.
This idea has been previously used 
in the definition of core vertices within 
\cp structures \cite{holme05, rombach14}, and paves the way for the definition
presented below.

\begin{defin} \label{defin-csp}
A \csp structure is a 3-partitioned connected graph with the following 
(non-empty) vertex classes:\vspace{-2mm}
\begin{itemize}
\item[{\rm (i)}] {\em core vertices}, with eccentricity not greater than two;\vspace{-1mm}
\item[{\rm (ii)}] {\em \s vertices}, adjacent (at least) to a pair of non-adjacent vertices from 
the other two classes; and \vspace{-1mm}
\item[{\rm (iii)}] {\em \p vertices}, with degree one.\vspace{-2mm}
\end{itemize}
Moreover, the graph is required not to have proper 
T-twin core vertices or proper F-twin \spe subgraphs.
\end{defin}


Here, semiperiphery vertices are simply required to act as intermediaries between (at least)
a core and a periphery, whereas for the latter we impose a minimal connection to the rest of the
network, in a way which implies in particular that 
periphery vertices are isolated from each other
(cf.\ item (ii) on page \pageref{peripheryvertices}).
Note that the requirements depicted for each class may be satisfied by vertices from other classes: 
e.g.\ a core may have degree one and/or connect a pair of (non-adjacent) \s and periphery vertices, whereas 
a \s or a \p vertex might well have eccentricity not greater than two. It is pretty clear, however,
that the requirements in items (ii) and (iii) are mutually exclusive.

It is worth emphasizing that this approach admits further extensions; on the one hand
we may consider the maximum core eccentricity (mce) and maximum \p degree (mpd) 
as parameters which in our present framework 
are fixed to the values two and one in (i) and (iii), respectively. 
Allowing these parameters to take on higher values may well lead to other structures of interest.
Additionally, in a setting with mce $\geq 3$ 
we might also define structures with more than three (ranked) classes, by distinguishing 
several semiperiphery layers defined by vertices which are adjacent to
vertex pairs coming from a higher-rank and a lower-rank class
(examples of networks with four classes can be found in \cite{deNooy, lazega}). These ideas
define tentative lines for future research.

The twin-free conditions stated at the end of Definition \ref{defin-csp}, supported
on the ideas discussed in Sections \ref{sec-falsetwins}
and \ref{sec-truetwins}, are the key element to reduce the
seemingly large number of CSP structures. As already indicated in 
the Introduction and in subsection \ref{subsec-cp}, 
the core
should be thought as a set of heavily interconnected vertices, amounting to a
fully connected set in idealized  cases; for this reason the 
true-twin notion for vertices is enough to reduce the eventual number
of core subgraphs within \csp structures. On the other hand, the F-twin concept for the \spe subgraph arises
as a natural extension of the false-twin notion for periphery vertices discussed in subsection \ref{subsec-cp},
allowing one to reduce the number of \spe subgraphs as well.
Note also that the the non-adjacency property stated in Corollary \ref{coro-nonadj} captures
the fact that twin \spe substructures to be reduced should be somehow independent, being related
only through the core vertices; in other words, if two (or more) semiperiphery vertices are adjacent then
it is natural to consider them as part of the same substructure. 


\subsection{Decomposition 
of CSP structures}
\label{subsec-decomposition}

Definition \ref{defin-csp} allows for an explicit description of
\csp structures, as detailed below. 

\begin{theor} \label{th-cspcharact}
Core-semiperiphery-periphery structures meeting Definition
\ref{defin-csp} admit the decomposition described in the sequel.\vspace{-3mm}
\begin{enumerate}
\item The core subgraph ${\cal C}$ is a join ${\cal C}_0 + {\cal C}_1$,
where\vspace{-1mm}
\begin{itemize}
\item ${\cal C}_0$ is a complete graph $K_{n_0}$; and
\item ${\cal C}_1$ is any graph of order $n_1$ 
without T-twin vertices.
\end{itemize}
\item The core-semiperiphery subgraph is a join ${\cal C} + {\cal S}$, where
${\cal C}$ has the form described above and
${\cal S}$ is any graph or order $n_s$ without F-twin subgraphs.
\item The periphery subgraph ${\cal P}$ is an empty graph
of order $n_p=n_0 + n_s$. Periphery vertices are leaves 
attached in a one-to-one
basis either to a vertex from ${\cal C}_0$ or 
from  ${\cal S}$.\vspace{-2mm}
\end{enumerate}
The orders $n_c=n_0+n_1$, $n_s$ and $n_p$ do not vanish, but
either $n_0$ or $n_1$ 
may do. 
\end{theor}

\noindent {\bf Proof.} Note in advance that the splitting of core vertices
in two groups ${\cal C}_0$ and ${\cal C}_1$ is defined from the fact that
those in ${\cal C}_0$ are connected to a periphery vertex
whereas those in ${\cal C}_1$ are not, as stated in item 3. 
In this regard, it is obvious that
periphery vertices are only connected either to a core (in ${\cal C}_0$) 
or to a
semiperiphery vertex because of 
the degree one condition stated
in item (iii) of Definition \ref{defin-csp}; notice that a 
single $K_2$ consisting of two
peripheries is ruled out by the requirement that the graph has at
least one core and one semiperiphery vertex.
Conversely, semiperiphery vertices are necessarily
connected to a single 
periphery (in addition to
cores and, possibly, other semiperipheries), since
two or more peripheries eventually connected to the 
same semiperiphery vertex would be
false twins. For the same reason, a core vertex in ${\cal C}_0$ 
is attached
to one periphery (again, in addition to connections to other
cores and to semiperipheries). These properties fully describe
the structure of the periphery subgraph ${\cal P}$ and will be used
throughout the rest of the proof.

Regarding the structure of the core subgraph,
${\cal C}_0$ is a complete graph (maybe the null
one $K_0$) and, moreover, it defines
a join with (i.e.\ it is fully connected to) 
${\cal C}_1$, if non-empty, because of
the eccentricity requirement for core vertices.
Indeed,
suppose there is a pair of non-adjacent core vertices, at least one of
which is adjacent to a periphery (i.e.\ at least
one of which is in ${\cal C}_0$): the distance of this periphery vertex
to the other core in that pair would be at least three, against the 
assumption that the maximum eccentricity of core vertices is two 
as stated in item (i) of Definition \ref{defin-csp}.

The core and the semiperiphery are fully connected as well. Again,
assuming the contrary, the distance between
such a core and the periphery vertex adjacent to that semiperiphery would
be greater than two, against the aforementioned eccentricity requirement.

It remains to show that the exclusion of twin structures in Definition
\ref{defin-csp} is equivalent to the absence of the corresponding
twin structures in the core or semiperiphery subgraph,
respectively, in the terms
stated in this Proposition. 
Regarding core vertices, 
note first that ${\cal C}_0$
may never include T-twins (meant in the full graph)
since the peripheries attached to these cores are adjacent only 
to one core and, therefore, these peripheries necessarily make a difference
in the neighborhoods of the corresponding cores; for the
same reason, cores in ${\cal C}_0$ and in ${\cal C}_1$ may never
be T-twins in the full graph. Additionally,
the absence of T-twins in ${\cal C}_1$ 
can be equivalently checked in the full graph or in the
core subgraph because of the fact that cores in ${\cal C}_1$
are not adjacent to any peripheries and, on the
contrary, fully connected to
both ${\cal C}_0$ and ${\cal S}$; this means
that the neighborhoods of two ${\cal C}_1$-cores in the full
graph differ if and only if 
these core vertices have different neighbors within ${\cal C}_1$. 


Concerning the equivalence between F-twin structures, let us first assume that
two subgraphs $H_1$ and $H_2$ within the \spe subgraph are F-twins in the full
graph, and let $\varphi$ denote the
corresponding isomorphism, so that (\ref{Ns}) holds
for all $u \in V_1=V(H_1)$. 
Let $\varphi_s$ stand for the restriction of this isomorphism
to $H_1 \cap {\cal S}$, and denote $V_{1s}=V_1 \cap V({\cal S})$,
$V_{2s}=V_2 \cap V({\cal S})$. From (\ref{Ns}) 
we get
\begin{equation*} 
({\cal N}(u)- V_1)\cap V({\cal S}) = ({\cal N}(\varphi(u))- V_2)\cap V({\cal S}),
\end{equation*}
an identity that can be recast as
\begin{equation} \label{Nsnn}
{\cal N}_s(u)- V_{1s} = {\cal N}_s(\varphi(u))- V_{2s}
\end{equation}
by making use of the property $(A-B) \cap C = A \cap C - B \cap C$ for 
arbitrary sets $A$, $B$, $C$ (here ${\cal N}_s(u)$ denotes ${\cal N}(u)
\cap V({\cal S})$).
By noting that (\ref{Nsnn}) holds
for all $u \in V_{1s}$ and that $\varphi(u)=\varphi_s(u)$ for vertices
in $V_{1s}$, it follows that $H_1 \cap {\cal S}$ and  
$H_2 \cap {\cal S}$ are F-twins as
subgraphs of ${\cal S}$ via the restricted isomorphism
$\varphi_s$, as we aimed to show.

Conversely, let $H_{1s}$ and $H_{2s}$ be F-twin structures as
subgraphs of ${\cal S}$, and denote by $\varphi_s$ the corresponding
isomorphism. Denote by $V_{1s}$ and $V_{2s}$ the vertex sets
of  $H_{1s}$ and $H_{2s}$, respectively. Let $H_1$ (resp.\ $H_2$) be the
subgraph induced 
 in the full graph 
by the vertices of $V_{1s}$ (resp.\ $V_{2s}$) and their adjacent
peripheries, and write as $V_1$ (resp.\ $V_2$) be the vertex set of 
$H_1$ (resp.\ $H_2$). Now, for every  $u \in V({\cal S})$ write as 
$p(u)$ the unique periphery vertex attached
to $u$ in the full graph and, conversely,
for every $u \in {\cal P}$ let $s(u)$ be
the unique semiperiphery vertex adjacent to $u$. With this notation 
we extend the isomorphism $\varphi_s$
to the whole of $H_1$ by setting
$$\varphi(u)= \begin{cases} 
\varphi_s(u) & \text{ if } u \in H_1 \cap {\cal S} \\
p(\varphi_s(s(u))) & \text{ if } u \in H_1 \cap {\cal P}.\vspace{1mm}
\end{cases}
$$
We claim that $\varphi$ makes $H_1$ and $H_2$ F-twin subgraphs in the 
full graph.
First, note that by construction (\ref{Nsnn}) is met for all
$u \in H_1 \cap {\cal S}$, and then  
\begin{equation} \label{Nsnnbis}
\big(V({\cal C}) \cup {\cal N}_s(u)\big)- V_{1s} = 
\big(V({\cal C}) \cup {\cal N}_s(\varphi(u))\big)- V_{2s}
\end{equation}
holds because $V({\cal C}) \cap V_{1s} = V({\cal C}) \cap V_{2s}
=\emptyset$; additionally, since 
$V_1 - V_{1s}$ and $V_2 - V_{2s}$ are in the periphery, 
we may rewrite (\ref{Nsnnbis}) as
\begin{equation} \label{Nsnnter}
\big(V({\cal C}) \cup {\cal N}_s(u)\big)- V_1 = 
\big(V({\cal C}) \cup {\cal N}_s(\varphi(u))\big)- V_2.
\end{equation}
Moreover, using the fact that $p(u) \in V_1, \ p(\varphi(u)) \in V_2$,
(\ref{Nsnnter}) yields
\begin{equation} \label{Nsnncuat}
\big(V({\cal C}) \cup {\cal N}_s(u) \cup \{p(u)\} \big)- V_1 = 
\big(V({\cal C}) \cup {\cal N}_s(\varphi(u)) \cup \{p(\varphi(u))\} \big)- V_2.
\end{equation}
In light of
the join structure proved above for ${\cal C}+{\cal S}$
we have
${\cal N}(u) = V({\cal C}) \cup {\cal N}_s(u) \cup \{p(u)\}$
and ${\cal N}(\varphi(u)) = V({\cal C}) \cup {\cal N}_s(\varphi(u)) 
\cup \{p(\varphi(u))\}$
for every $u \in H_1 \cap {\cal S}$, so that (\ref{Nsnncuat}) is equivalent
to (\ref{Ns}). 

It remains to show that (\ref{Ns}) 
also holds for $u \in H_1 \cap {\cal P}$, but this
is a much simpler check. Indeed, we have ${\cal N}(u) = \{s(u)\}$ and, by construction, 
$s(u) \in V_1$, so that the left-hand side of (\ref{Ns}) 
is ${\cal N}(u) -V_1 = \emptyset$.
Analogously, $\varphi(u) = p(\varphi_s(s(u)))$ and therefore 
${\cal N}(\varphi(u)) = \{\varphi_s(s(u))\}$; again, $\varphi_s(s(u)) \in V_2$ and 
the right-hand side of (\ref{Ns}) 
also verifies ${\cal N}(\varphi(u))-V_2 = \emptyset.$ This
means that (\ref{Ns}) 
holds trivially if $u \in H_1 \cap {\cal P}$ and this, together
with the remarks in the previous paragraph, shows that $H_1$ and $H_2$ as constructed above
are F-twins in the full graph.





Note finally 
that, apart from the twin-free requirements above, 
both ${\cal C}_1$ and ${\cal S}$ admit any 
topology since no additional restrictions emanate from Definition
\ref{defin-csp}. This completes the proof of Theorem \ref{th-cspcharact}.

\hfill $\Box$


\subsection{Enumeration of CSP structures} 
\label{subsec-enumeration}

Theorem \ref{th-cspcharact} above essentially 
reduces the enumeration problem for CSP structures
to a combination of a subgraph 
${\cal C}_1$ within the core
displaying no
true twin vertices, and a \s subgraph ${\cal S}$ without any kind of
F-twins, 
with the eventual addition (join) of a complete graph ${\cal C}_0$
with its corresponding peripheries attached.
In this 
problem one is faced with two different sub-problems of independent
mathematical interest: enumerating 
graphs without true twin vertices on the one
hand, and graphs without F-twin subgraphs on the other.
 We let $t_n$ and $s_n$ be the numbers of graphs
on $n$ vertices without true twin vertices and without F-twin
subgraphs, respectively. 
It is worth mentioning 
that, in light of Theorem \ref{th-duality}, these two
numbers coincide with those of graphs without false twin vertices 
and graphs without T-twin subgraphs, although we will not make use
of this except for the obvious remark that $s_n \leq t_n$.
Related enumeration problems are finding the numbers 
of graphs without any type of twin vertices (that is, without
either true or false twin vertices) and without either T-twin or F-twin
subgraphs.

The number
of \csp structures  can be computed in arbitrary order 
($\geq 3$) in terms of
the quantities $t_n$ and $s_n$ defined above. 
We will do so by splitting the computation
in two parts.
First we compute the number $x_n$ of \csp structures of order $n$ in which
all periphery vertices are adjacent to the semiperiphery:
this corresponds
to the case $n_0 = 0$ (or ${\cal C}_0=K_0$) in the notation of Theorem \ref{th-cspcharact}. 
Later on 
we will add a number $y_n$ of structures with $n_0 > 0$ to get the total
number $z_n=x_n+y_n$ of CSP structures on $n$ vertices.

In order to compute $x_n$, by means of Theorem \ref{th-cspcharact}
the number of joins ${\cal C}_1 + {\cal S}$ is easily 
seen to be given by all combinations of $t_{n_c}$ 
core subgraphs on $n_c$ vertices without true twins
and $s_{n_s}$ semiperiphery subgraphs on $n_s$ vertices without F-twin subgraphs.
Using the fact that in this setting $n_s=n_p$ and then
$n=n_c + 2n_s$,
some easy computations yield
\begin{equation}\label{csp-nosat}
x_n  = \begin{cases}
\dsp \sum_{k=1}^{\frac{n-1}{2}}t_{2k-1}s_{\frac{n+1}{2}-k} & \text{ if } n \text { is odd} \\
\vspace{-2mm} \\
\dsp \sum_{k=1}^{\frac{n-2}{2}}t_{2k}s_{\frac{n}{2}-k} & \text{ if } n \text { is even,} \\
\end{cases}
\end{equation}
for $n \geq 3$.

On the other hand, we can compute $y_n$ in a recursive manner, just using the remark
that all structures with $n_0 > 0$ can be obtained from a lower order structure 
just joining (the core vertex of) a core-periphery pair to the cores and semiperipheries of
this lower order structure.
This leads to 
\begin{equation}\label{csp-sat}
y_n  = \begin{cases}
z_{n-2} & \text{ if } n \text { is odd} \\
z_{n-2} + s_{\frac{n}{2}-1}& \text{ if } n \text { is even,} \\
\end{cases}
\end{equation}
again for $n \geq 3$. The additional term $s_{\frac{n}{2}-1}$
for even $n$ captures the structures with only one 
core which belongs to ${\cal C}_0$.
Note that we make recursive 
use of the total number $z_n = x_n + y_n$ of \csp structures,
setting $z_1 = z_2 = 0$ for consistency. 

Equations (\ref{csp-nosat}) and (\ref{csp-sat}) together 
define recursively the total number of \csp structures on $n$ vertices, which
(omitting details for the sake of brevity) read,
%
in terms of the numbers $n$ (total number of vertices) 
and $n_c$ (number of core vertices), as 
\begin{equation}\label{csp-total-tris}
z_{n,n_c}  = \begin{cases}
\dsp \sum_{k=0}^{\min \{ E(\frac{n_c-1}{2}),\frac{n-n_c-3}{2} \}}  
t_{n_c-2k-1} s_{\frac{n-n_c-1}{2}-k}& \text{ if } n -n_c \text { is odd} \\
\vspace{-2mm} \\
\dsp \sum_{k=0}^{\min \{ E(\frac{n_c}{2}),\frac{n-n_c}{2}-1 \}}
t_{n_c-2k} s_{\frac{n-n_c}{2}-k} & \text{ if } n-n_c \text { is even.} \\
\end{cases}
\end{equation}
Finally, $z_n$ is the sum of the above values of $z_{n,n_c}$ for $n_c=1 \ldots n-2$.


\begin{table}[ht]  \label{table-notruetwins}
\caption{Number of graphs without (true) twin vertices}
\begin{center}
\begin{tabular}{cc}
\toprule
\ \ Order ($n$) \ \ & \ \ \#Graphs without (true) twin vertices ($t_n$) \ \
\\ \midrule
1 & {\bf 1} \\ \midrule
2 & {\bf 1} \\ \midrule
3 & {\bf 2} \\ \midrule
4 & {\bf 5} \\ \midrule
5 & {\bf 16} \\\midrule
6 & {\bf 78} \\ 
\bottomrule
\end{tabular}
\end{center}
\end{table}

In the sequel we use the above derived formulas to compute the number of 
core-semiperi\-phery-periphery structures in low order (up to $n= 8$), in terms of the 
previously defined quantities $t_n$ and $s_n$. To the knowledge of the author,
the number $t_n$ of graphs without true twin vertices 
(or without false twins vertices) is not
known in general; however, computationally this is a very simple task
in low order and for later use we depict the numbers $t_n$ up to $n=6$ 
in
Table 1. 


\begin{table}[ht]  \label{table-cspstructures}
\caption{Number of \csp structures}
\begin{center}
\begin{tabular}{cc}
\toprule
\ \ Order ($n$) \ \ & \ \ Number of CSP structures ($z_n$) \ \
\\ \midrule
3 & {\bf 1} \\ \midrule
4 & {\bf 2} \\ \midrule
5 & {\bf 4} \\\midrule
6 & {\bf 9} \\\midrule
7 & {\bf 24} \\\midrule
8 & {\bf 96} \\ 
\bottomrule
\end{tabular}
\end{center}
\end{table}

The computation of $s_n$ (that is, the number of graphs on $n$ vertices
without any kind of F-twin subgraphs)
is more involved even from a computational point of view. Nevertheless, it is very
easy to check that the lowest order structure involving F-twin subgraphs
with order greater than one 
is $K_2 \cup K_2$; this obviously implies that $s_n =t_n$ for $n \leq 3$.
Additionally, one can easily see that only 
the subindices $i=1, 2, 3$ for $s_i$
are involved
in the computation of the number of CSP structures up to order eight. 
Using these remarks, the numbers $z_n$ up to $n=8$ 
are given in Table 2. 

\subsection{CSP structures in low order}
\label{subsec-examples}

The 
\csp structures in order up to 6 are displayed in Figures 
3 and 4. 
Core, semiperiphery and periphery vertices are painted black, 
grey
and white, respectively. Worth commenting are the facts that with $n=3$ one gets
the expected ``elementary'' CSP structure, and that
one of the two cases with $n=4$ arises from the 
addition of a periphery vertex connected to a (say) ${\cal C}_0$ core
vertex; a structure with two cores is already displayed 
in order four. Note also that up to three and four 
cores are displayed with $n=5$ and $n=6$.


\begin{figure}[ht] \label{fig-cspn3}
\vspace{6mm}
\mbox{} \hspace{-1mm}
\parbox{0.5in}{
\epsfig{figure=fig_csp3.eps, width=0.22\textwidth}
}
\hspace{38mm} 
\parbox{0.5in}{
\epsfig{figure=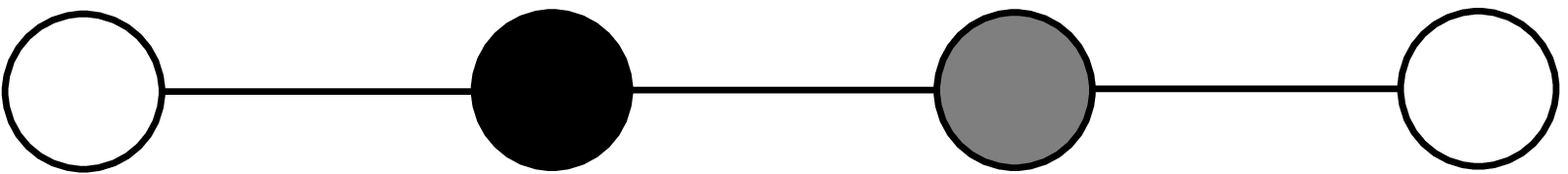, width=0.3\textwidth}
}
\hspace{53mm}
\parbox{0.5in}{
\epsfig{figure=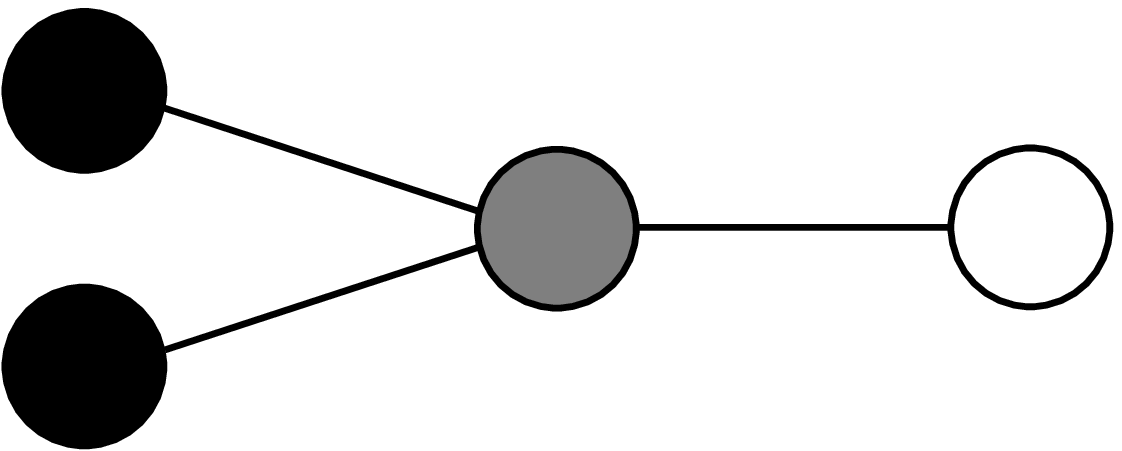, width=0.22\textwidth}
}\\ \vspace{16mm}\\
\mbox{} \hspace{23mm}
\parbox{0.5in}{
\epsfig{figure=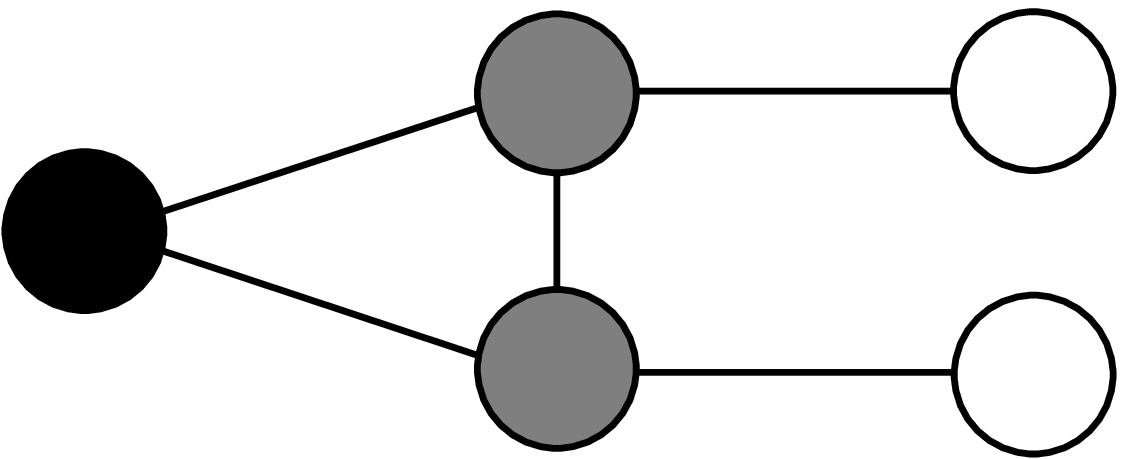, width=0.22\textwidth}
}
\hspace{45mm} 
\parbox{0.5in}{
\epsfig{figure=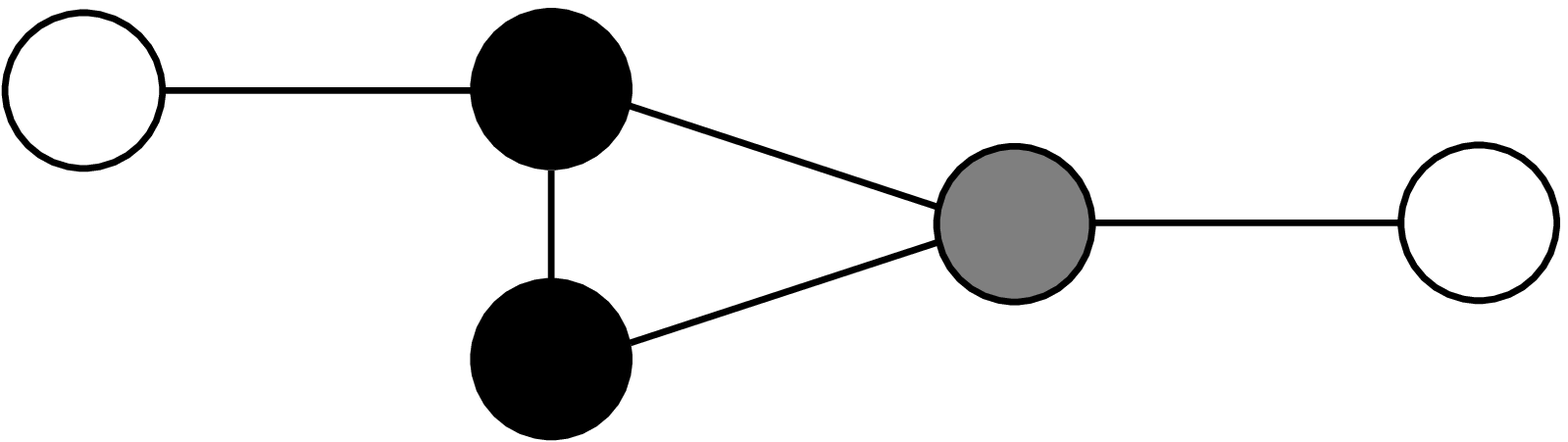, width=0.3\textwidth}
}
\vspace{14mm}\\
\mbox{} \hspace{23mm}
\parbox{0.5in}{
\epsfig{figure=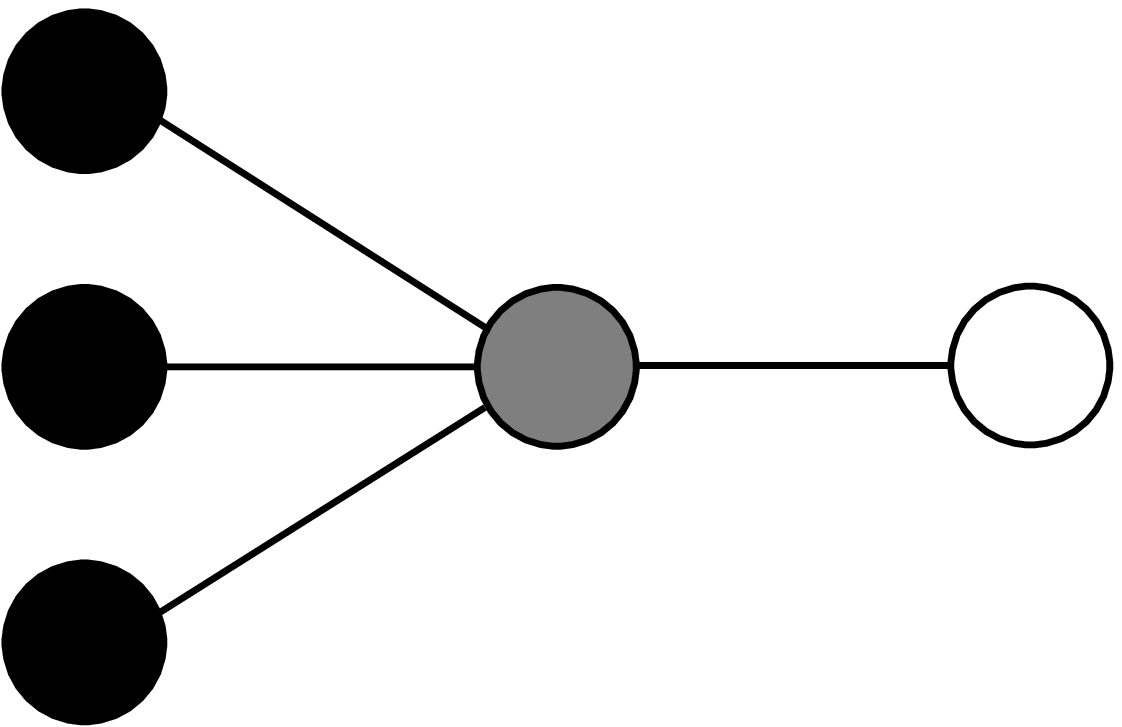, width=0.22\textwidth}
}
\hspace{60.7mm}
\parbox{0.5in}{
\epsfig{figure=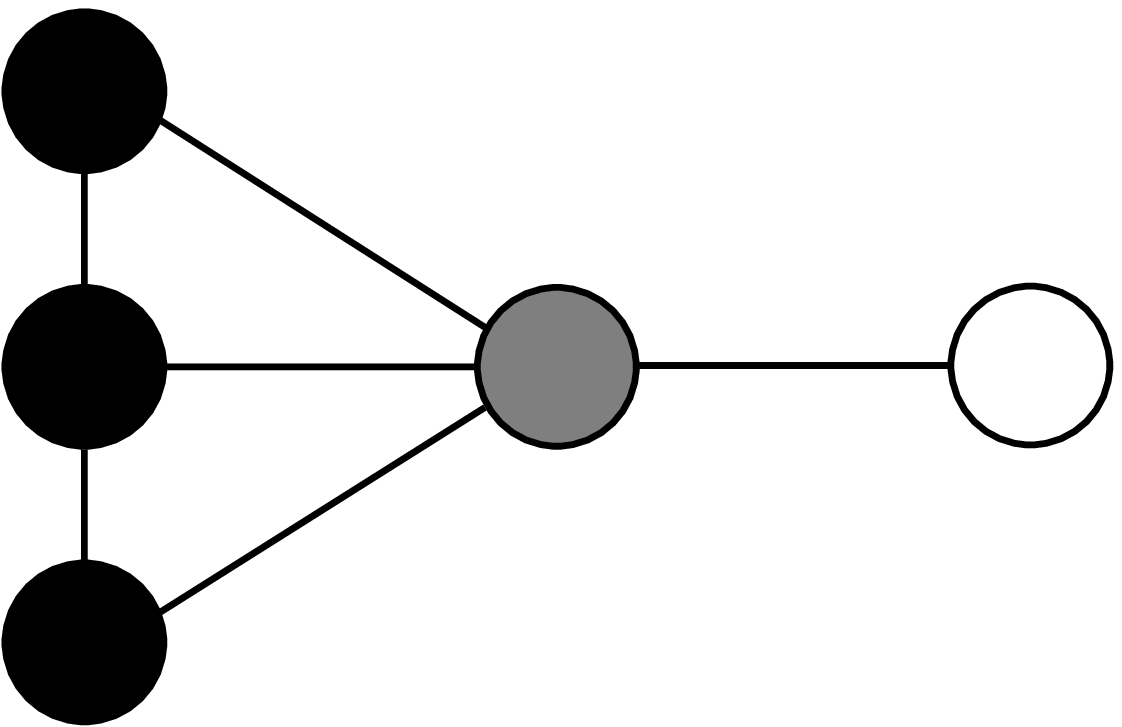, width=0.22\textwidth}
}
\vspace{6mm}
\caption{CSP structures up to order five} 
\end{figure}
\vspace{-3mm}
\begin{figure}[h]\label{fig-cspn6} 
\vspace{6mm}
\mbox{} \hspace{-8mm}
\parbox{0.5in}{
\epsfig{figure=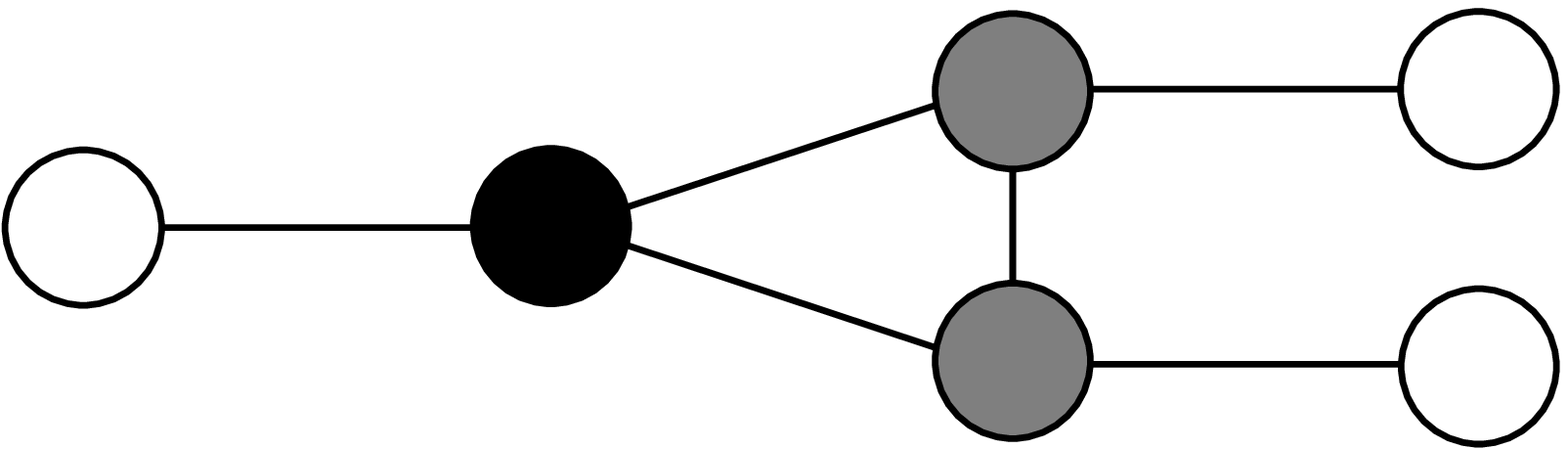, width=0.3\textwidth}
}
\hspace{53mm}
\parbox{0.5in}{
\epsfig{figure=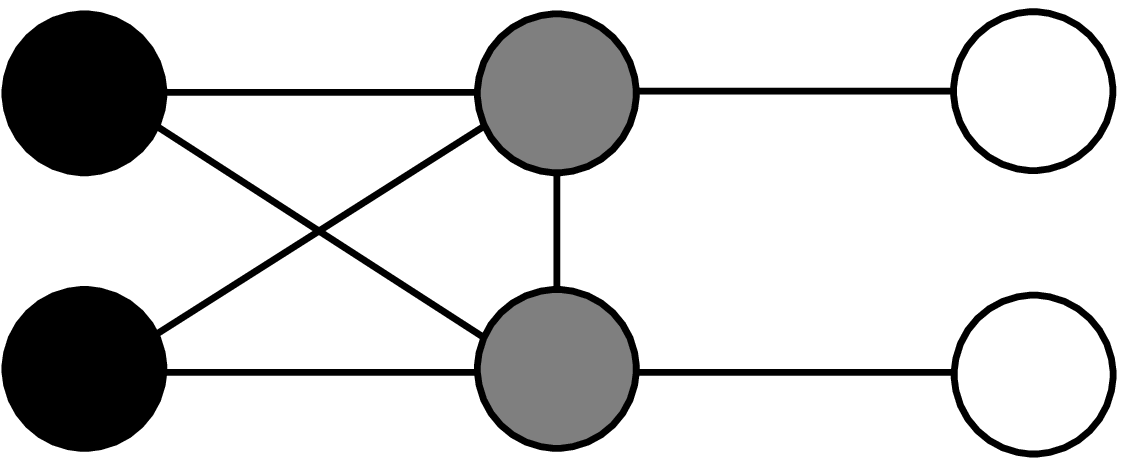, width=0.22\textwidth}
}
\hspace{40mm}
\parbox{0.5in}{
\epsfig{figure=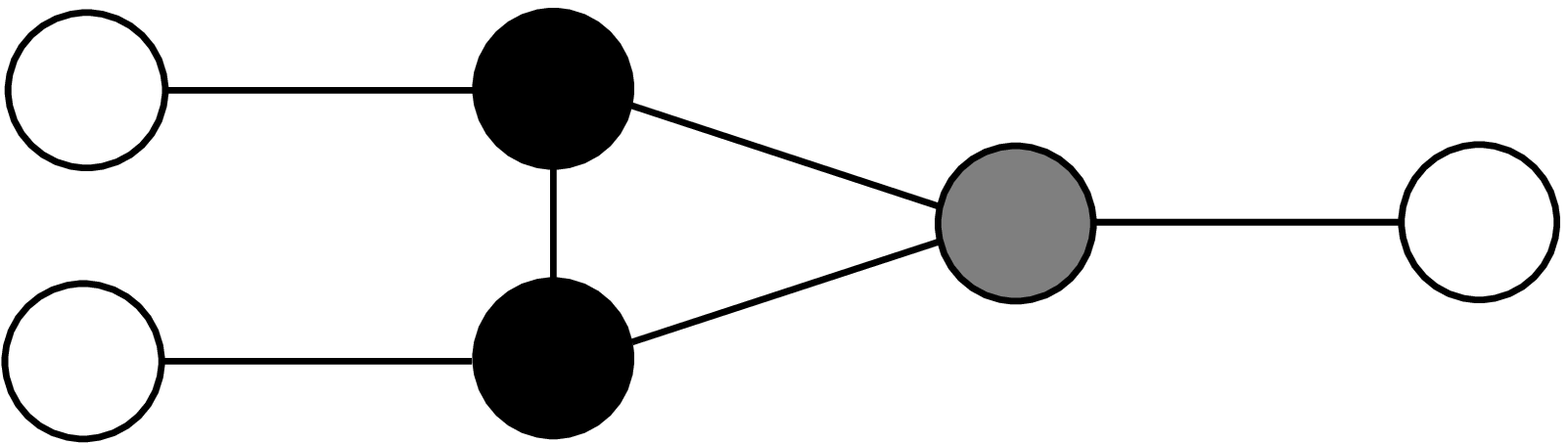, width=0.3\textwidth}
}\vspace{12mm}\\
\mbox{} \hspace{-8mm}
\parbox{0.5in}{
\epsfig{figure=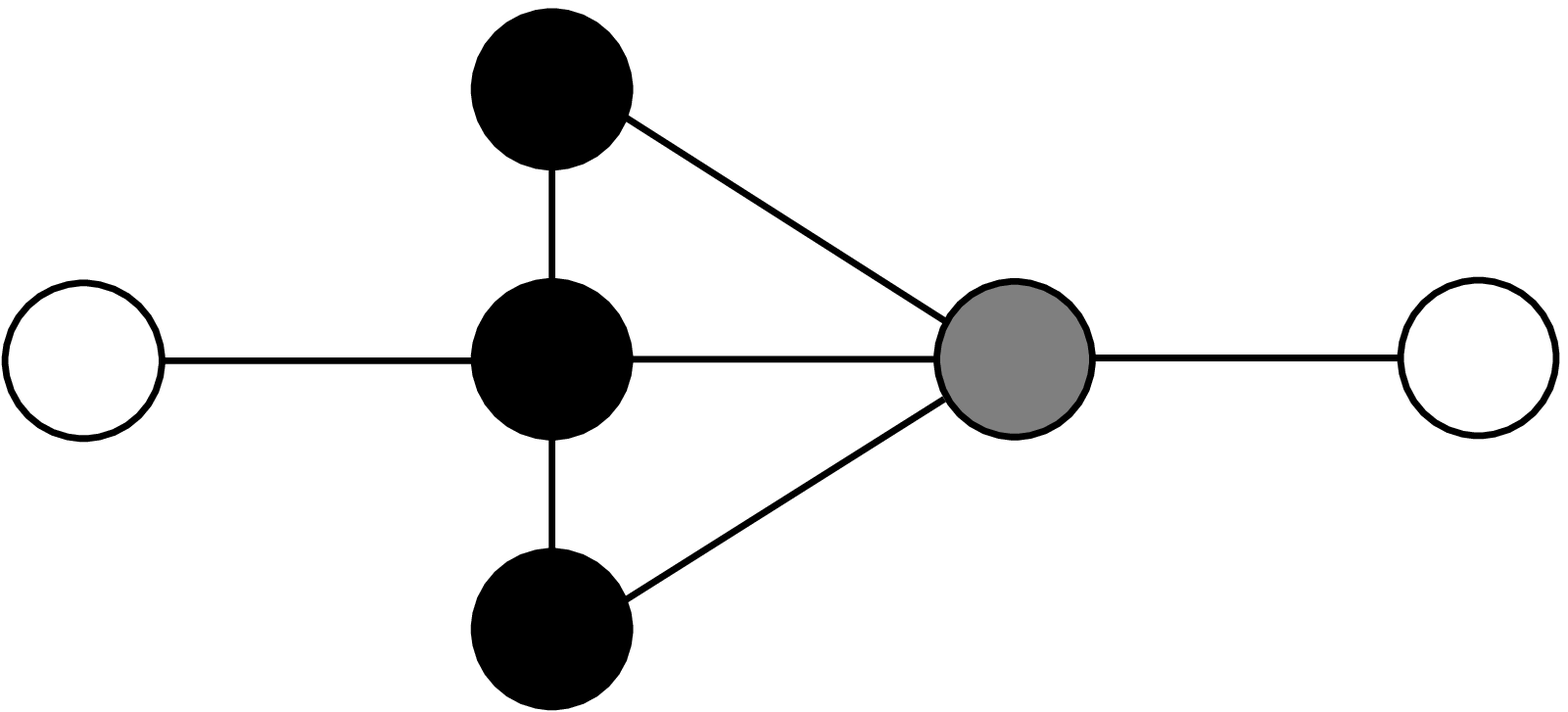, width=0.3\textwidth}
}
\hspace{47mm}
\parbox{0.5in}{
\epsfig{figure=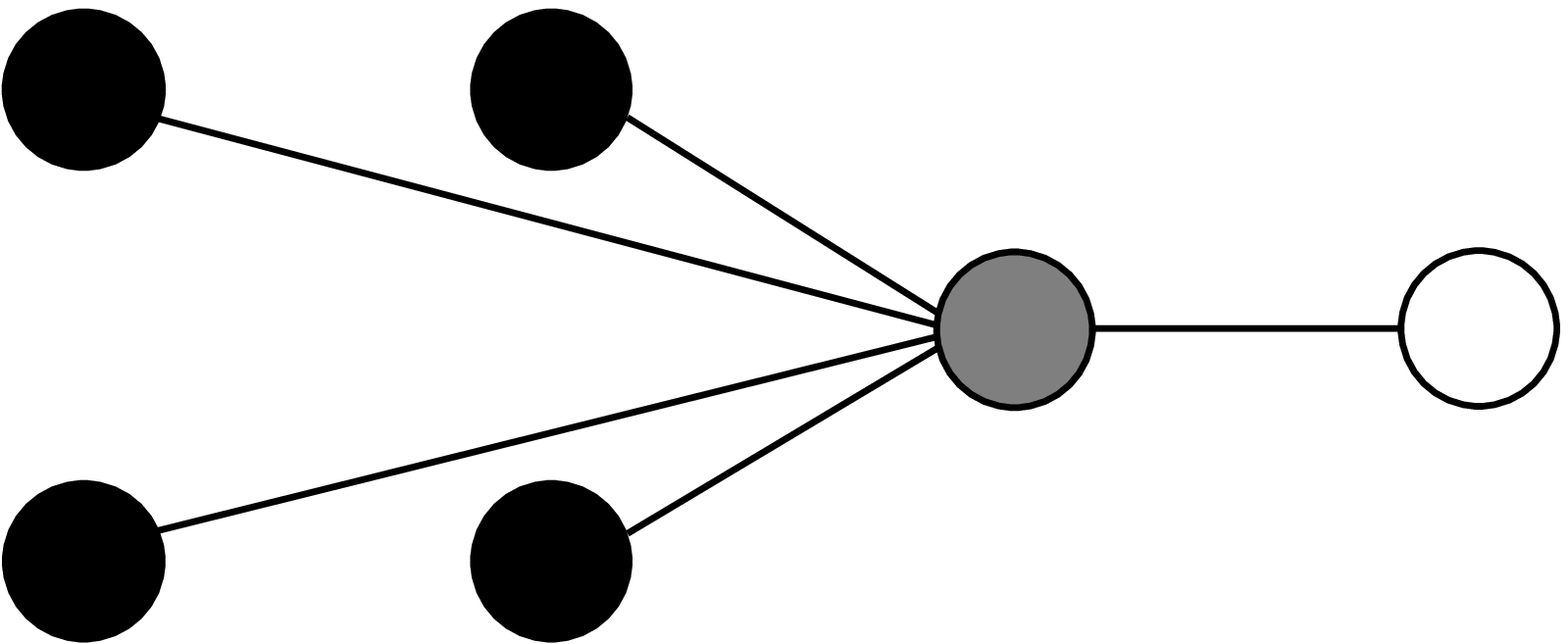, width=0.3\textwidth}
}
\hspace{47mm}
\parbox{0.5in}{
\epsfig{figure=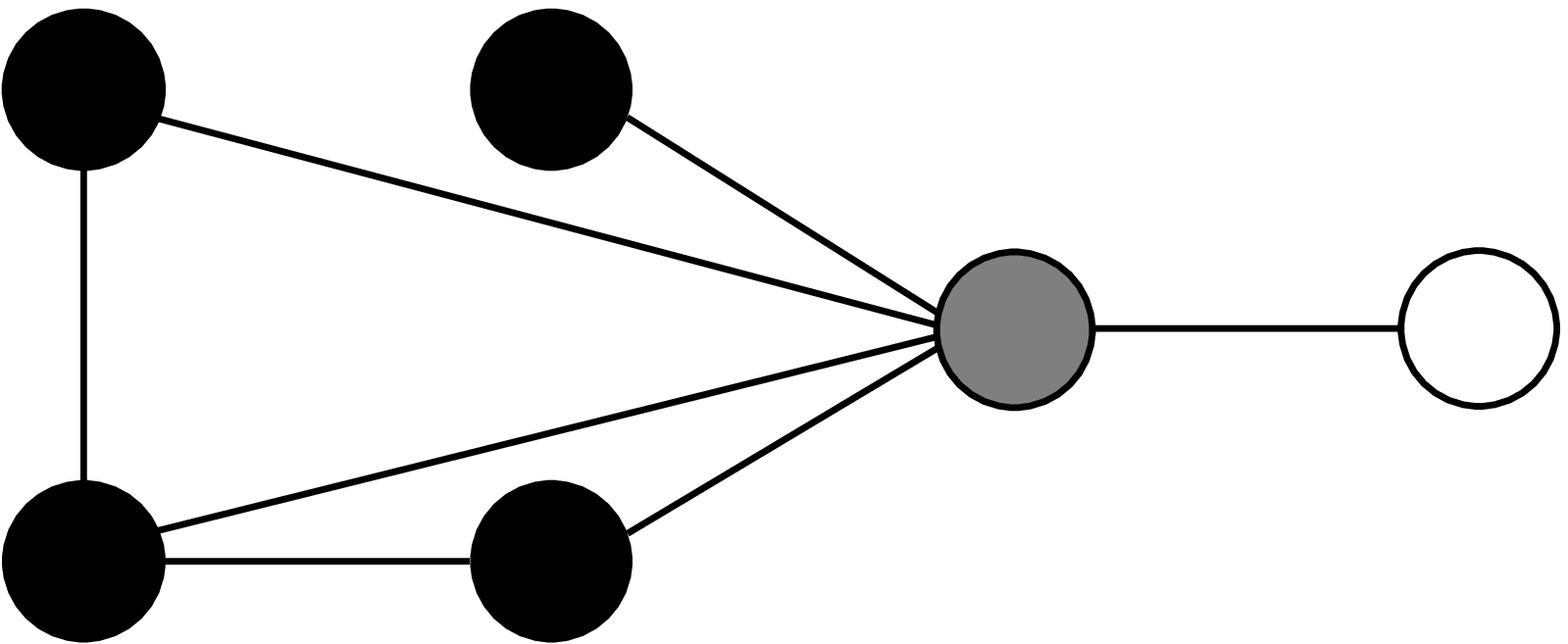, width=0.3\textwidth}
}\vspace{14mm}\\
\mbox{} \hspace{-8mm}
\parbox{0.5in}{
\epsfig{figure=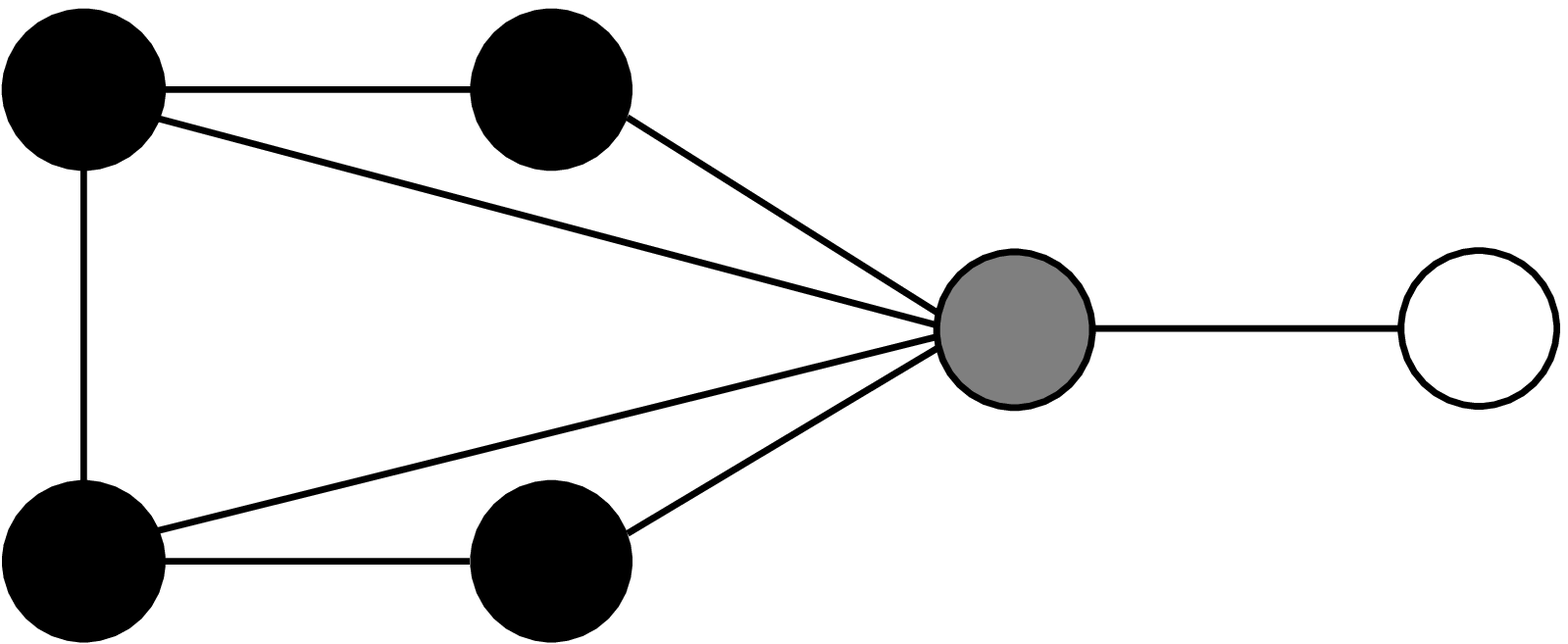, width=0.3\textwidth}
}
\hspace{47mm}
\parbox{0.5in}{
\epsfig{figure=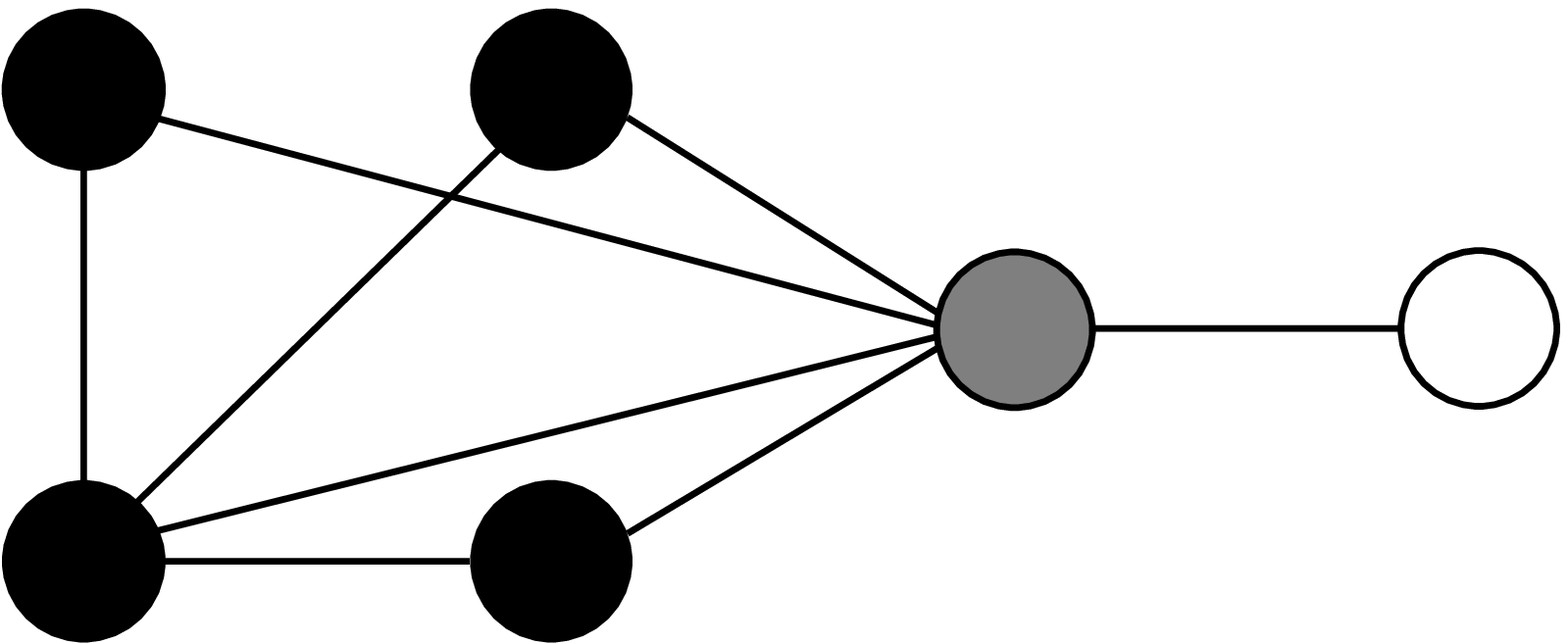, width=0.3\textwidth}
}
\hspace{47mm}
\parbox{0.5in}{
\epsfig{figure=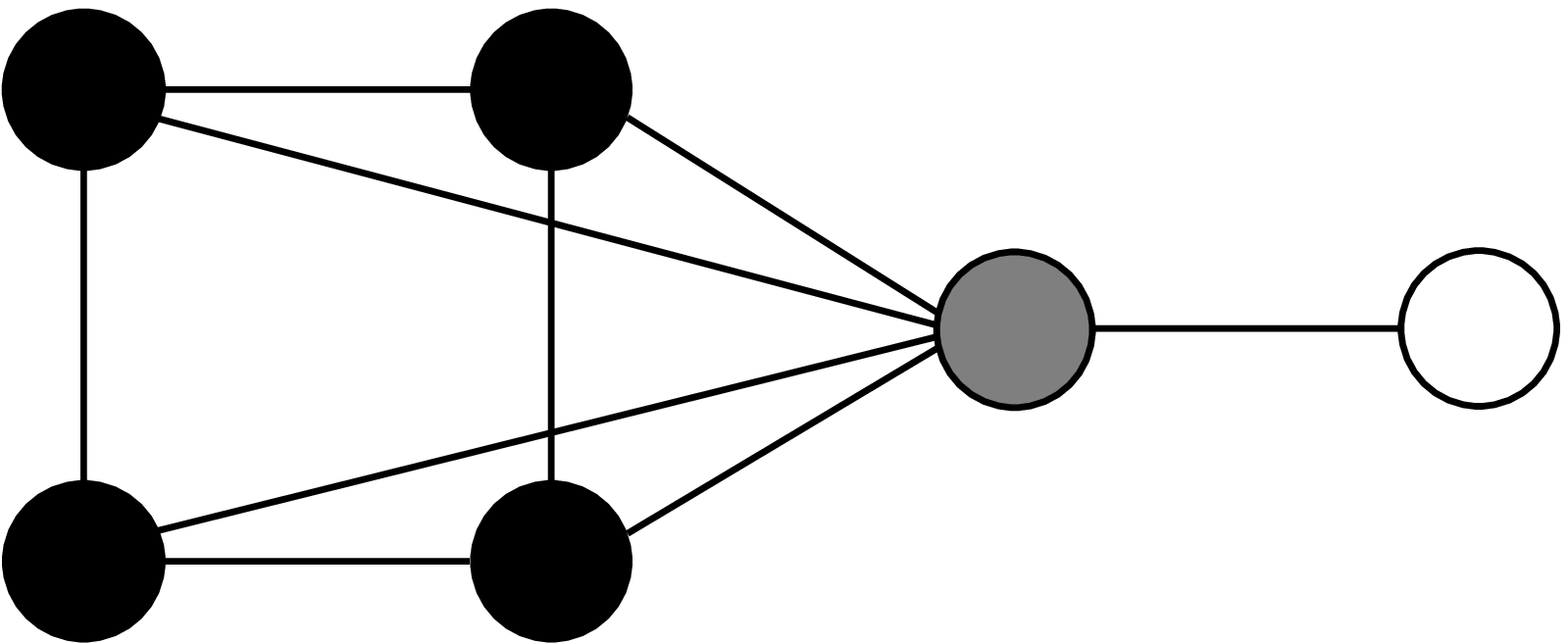, width=0.3\textwidth}
}
\vspace{6mm}
\caption{CSP structures in order six} 
\end{figure}


\section{CSP structure within the Asia-Africa-Oceania subnetwork of 1994 metal manufactures
trade}
\label{sec-ex}

The approach developed in previous sections
provides a formal definition and a criterion for the 
systematic classification of \csp structures in networks. 
In order to identify such structures in real problems, we need
to develop additional results based on positional analyses 
allowing one to assign systematically 
vertices to clusters and to evaluate the extent to which 
the quotient network fits a 
CSP structure. This task, in its broad 
generality, exceeds the scope of the present paper and will be the object 
of future research. However, we discuss below a roadmap for 
this research by examining a given subnetwork
of the network of miscellaneous imports of metal manufactures between 
80 countries in 1994. 
These data, coming from world trade statistics, have been 
previously addressed in \cite{deNooy} along the lines
discussed in the original work of Wallerstein \cite{wallerstein74}.
This data set is freely available on the web (cf.\ \cite{deNooy}).

Since the results in this section have illustrative purposes and
in order to simplify the discussion we restrict the attention
to a subnetwork of the abovementioned network, namely the
one defined by the
countries from Asia, Africa and Oceania for which data are available
in the original dataset.
Note that the large amount of exports of high-technology products from East Asian countries
makes 
this analysis relevant, looking in particular for their relation patterns 
with developing and least-developed countries from Africa, Oceania
and other regions of Asia.
In our model, every edge in the network is weighted with
the total amount of trade between the two countries (that is, we add
imports and exports). To reduce dimensionality
we remove edges in which this amount does not reach 10M (10 million)  USD 
or links involving countries whose total amount of trade does not reach 
25M USD; note that these quantities barely represent a few parts
per thousand of the total amount of trade in this network which is over 8 billion USD.
Exceptions are made when such a removal renders the network disconnected:
for the involved countries we then retain the edge displaying the highest amount of
trade with any of their commercial neighbors. This yields a connected network with
29 nodes and 69 edges (data are displayed on the Appendix).

In order to examine the presence  of CSP structures
in this network, as well as the eventual reduction 
of twin substructures, we use
two different
criteria to cluster vertices. The first one is very elementary and just
uses a threshold in the volume of trade between pairs of countries: 
we use 
this basic approach to provide simple examples of 
CSP structures and twin subgraphs. The second criterion
is more elaborate: in order to identify clusters we combine the amount 
of trade between countries, as above,
with a dissimilarity measure capturing
similar relation patterns. This will result in a refinement of the 
CSP structures which arise under the first clustering
criterion. Details are given 
below.

As indicated above, let us first cluster the different countries 
using the connected components of the graph which results
from removing edges below a given trade threshold. 
Let us for instance consider pairs of countries exchanging
at least 75M USD. This yields a main cluster defined by 11 countries,
namely China, Hong Kong, Japan, Thailand, Korea (to be referred in 
the sequel as East Asian countries), together with Malaysia, Singapore,
Indonesia, the Philippines 
(Southeast Asia), and Australia and New Zealand (both countries being jointly
referred to as Australasia). This
cluster comprises more than 7.7 billion 
USD trade, that is,
more than 95\% 
of the total amount of trade in the network.
None of the remaining countries reaches the above threshold with any neighbor,
so that each one of the 
other clusters is identified with a single country.

With this clustering,
the quotient graph displays 3 countries (Algeria, South Africa, and India)
which are adjacent to the main cluster and to 5 countries with degree one
(Tunisia (Algeria),  Israel, Mauritius, Reunion (South Africa),
and Oman (India), respectively).
There are 10 countries 
with degree one which are adjacent to the main cluster (Pakistan,
Bangladesh, Egypt, Jordan, Kuwait, Morocco, Madagascar, 
Seychelles, Sri Lanka and Fiji). This quotient network is displayed in
Figure 
5(a); we explicitly label the vertices
corresponding to Algeria, South Africa and India for better clarity.

\begin{figure} \label{fig-AAOceania1}
\vspace{6mm}
\hspace{-10mm}
\parbox{0.5in}{
\epsfig{figure=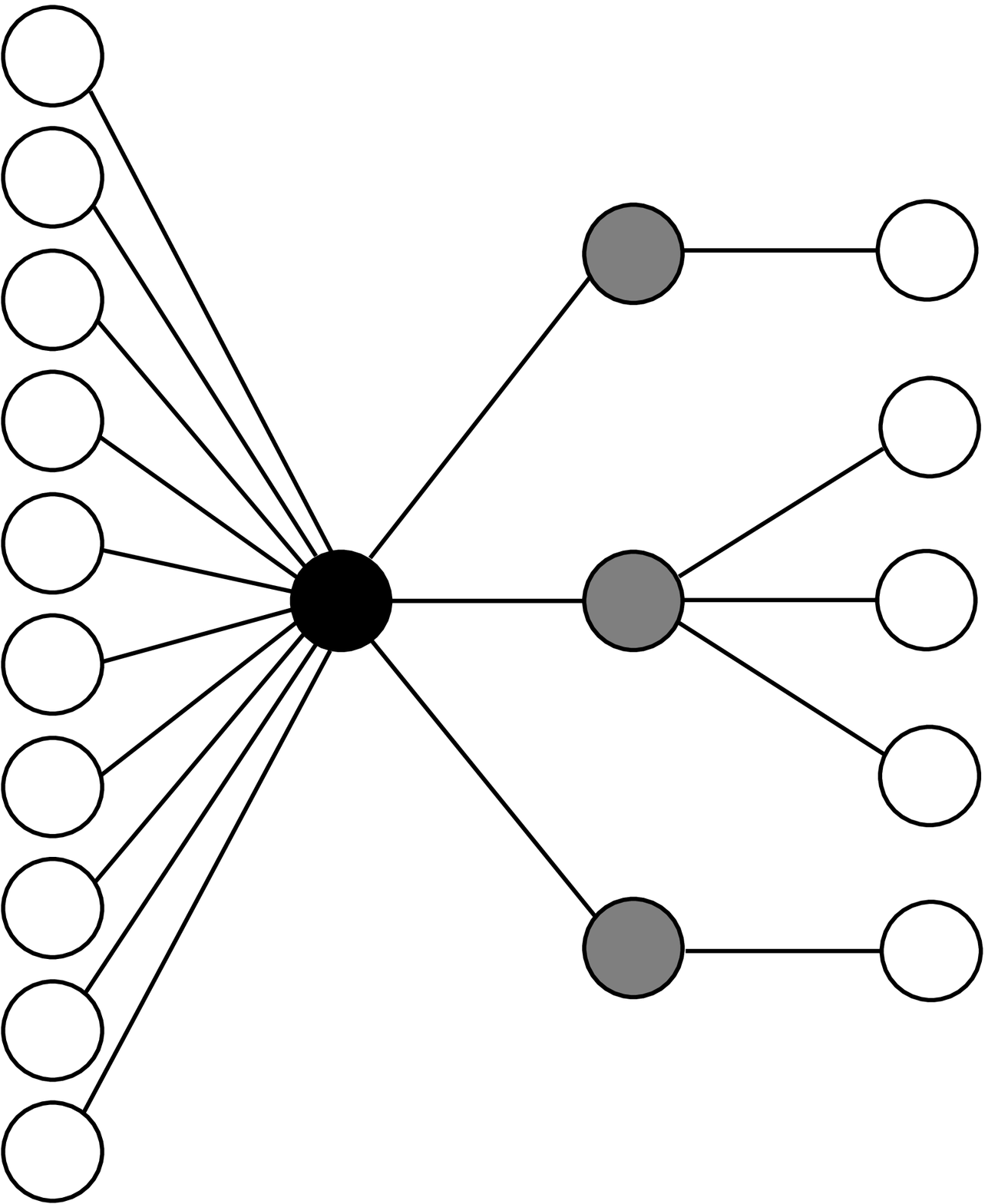, width=0.32\textwidth}
\put(-68,157){\footnotesize\sf Algeria}
\put(-73,105){\footnotesize\sf S.\hspace{-0.5mm} Africa}
\put(-59,51){\footnotesize\sf India}
}
\hspace{50mm} 
\parbox{0.5in}{
\epsfig{figure=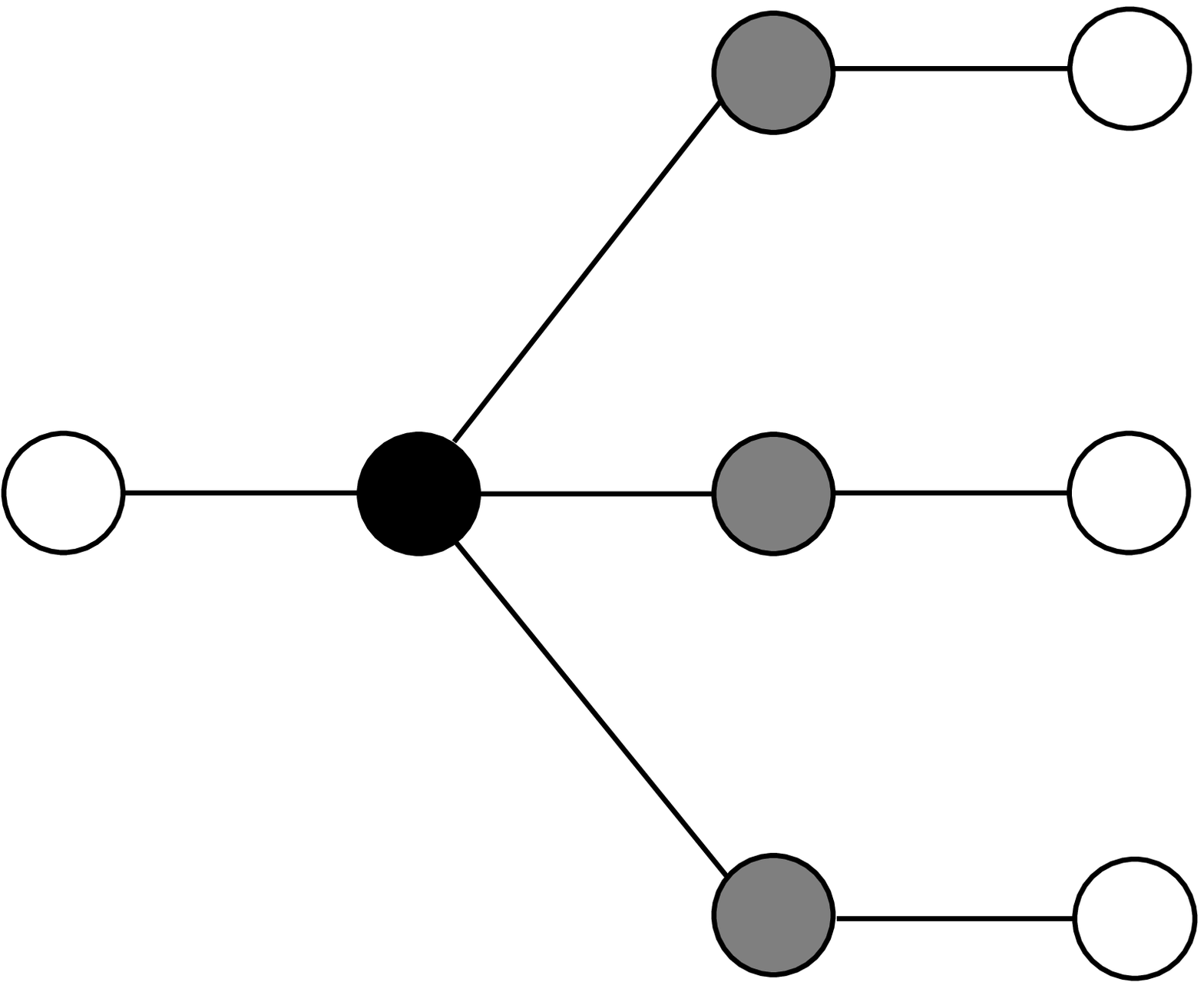, width=0.32\textwidth}
\put(-68,126){\footnotesize\sf Algeria}
\put(-73,73){\footnotesize\sf S.\hspace{-0.5mm} Africa}
\put(-57,20){\footnotesize\sf India}
\put(-24,126){\footnotesize\sf Tunisia}
\put(-20,73){\footnotesize\sf Israel}
\put(-21,20){\footnotesize\sf Oman}
\put(-154,74){\footnotesize\sf Pakistan}
\put(-111,74){\footnotesize\sf Core}
}
\hspace{50mm}
\parbox{0.5in}{
\epsfig{figure=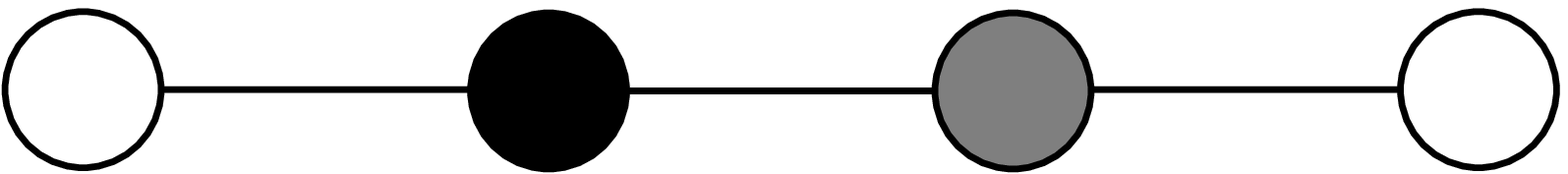, width=0.32\textwidth}
\put(-72,21.5){\footnotesize\sf S.\hspace{-0.5mm} Africa}
\put(-20,21.5){\footnotesize\sf Israel}
\put(-108,21.5){\footnotesize\sf Core}
\put(-155,21.5){\footnotesize\sf Pakistan}
}
\caption{Clustering with threshold 75M USD: (a) CSP network; (b) Reduction of false twin vertices; (c) CSP structure.}
\end{figure}


This quotient graph admits a classification of all the
clusters either as a core, semiperiphery or periphery, according
to the criteria given in Definition \ref{defin-csp}. 
The core is composed of the East and Southeast Asian countries
together with 
Australia and New Zealand, 
whereas
the semiperiphery is composed
of three countries (Algeria, South Africa, and India), and the fifteen 
countries 
with degree one define the periphery.
Among the latter,
the three ones adjacent to South Africa are false twins 
(we use Israel as their representative)
and,
analogously, the ten countries with degree one attached to the
core are false twins as well (with Pakistan as the
representative of this class). 
After identifying false twin vertices, the resulting graph is displayed
in Figure 5(b).
In turn, this figure clearly displays three subgraphs which are 
F-twins, namely,
the \spe pairs defined by Algeria and Tunisia, South Africa and Israel, and
India and Oman, respectively.  After identifying these three subgraphs
(with the pair South Africa-Israel being chosen as the representative of this relation
pattern),
the resulting CSP structure is depicted in 
Figure 5(c)
(it has four vertices and can be also found 
in Figure 3).
We emphasize that the F-twin notion makes it possible to
capture the elementary pattern displayed by the 
three \spe pairs mentioned above.


Another pattern arises if we raise 
the threshold to cluster countries
say to 125M USD. 
Since now neither Australia nor New Zealand trades such an amount
with any Asian country, but they do with each other, they
turn to define a cluster  by themselves (Australasia in the sequel),
independently of the 
East and Southeast Asian countries which are still joined
together into a big cluster, trading more than 7 billion USD. 
The latter still
meets the requirement defining a core in Definition \ref{defin-csp},
but the Australasian cluster does not, since it does not satisfy the 
eccentricity-two criterion (e.g.\ its distance to Israel is three).
Australasia may by contrast be classified as
a semiperiphery: note that Fiji is now attached to the Australasian cluster. 
The new quotient graph is displayed
in Figure 6(a).
As before,
we depict in Figure 
6(b) and (c), respectively, 
the network
without false twin vertices and the CSP
structure which finally results from removing
F-twin structures (now only the Algeria-Tunisia and South Africa-Israel
pairs).

\begin{figure} \label{fig-AAOceania2}
\vspace{6mm}
\hspace{-10mm}
\parbox{0.5in}{
\epsfig{figure=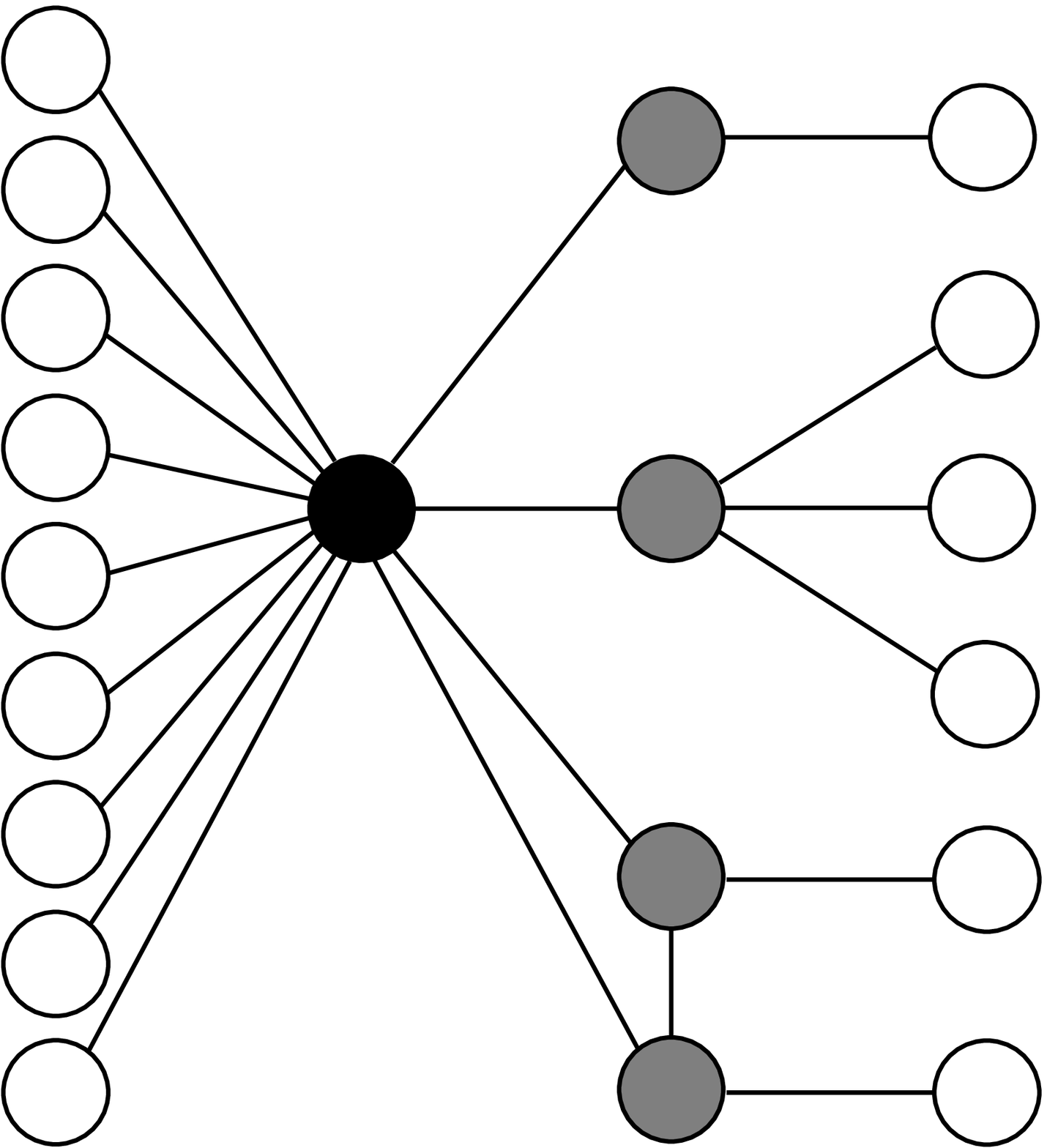, width=0.32\textwidth}
\put(-68,158){\footnotesize\sf Algeria}
\put(-74,105.5){\footnotesize\sf S.\hspace{-0.5mm} Africa}
\put(-60,51){\footnotesize\sf India}
\put(-78,-10){\footnotesize\sf Australasia}
}
\hspace{50mm} 
\parbox{0.5in}{
\epsfig{figure=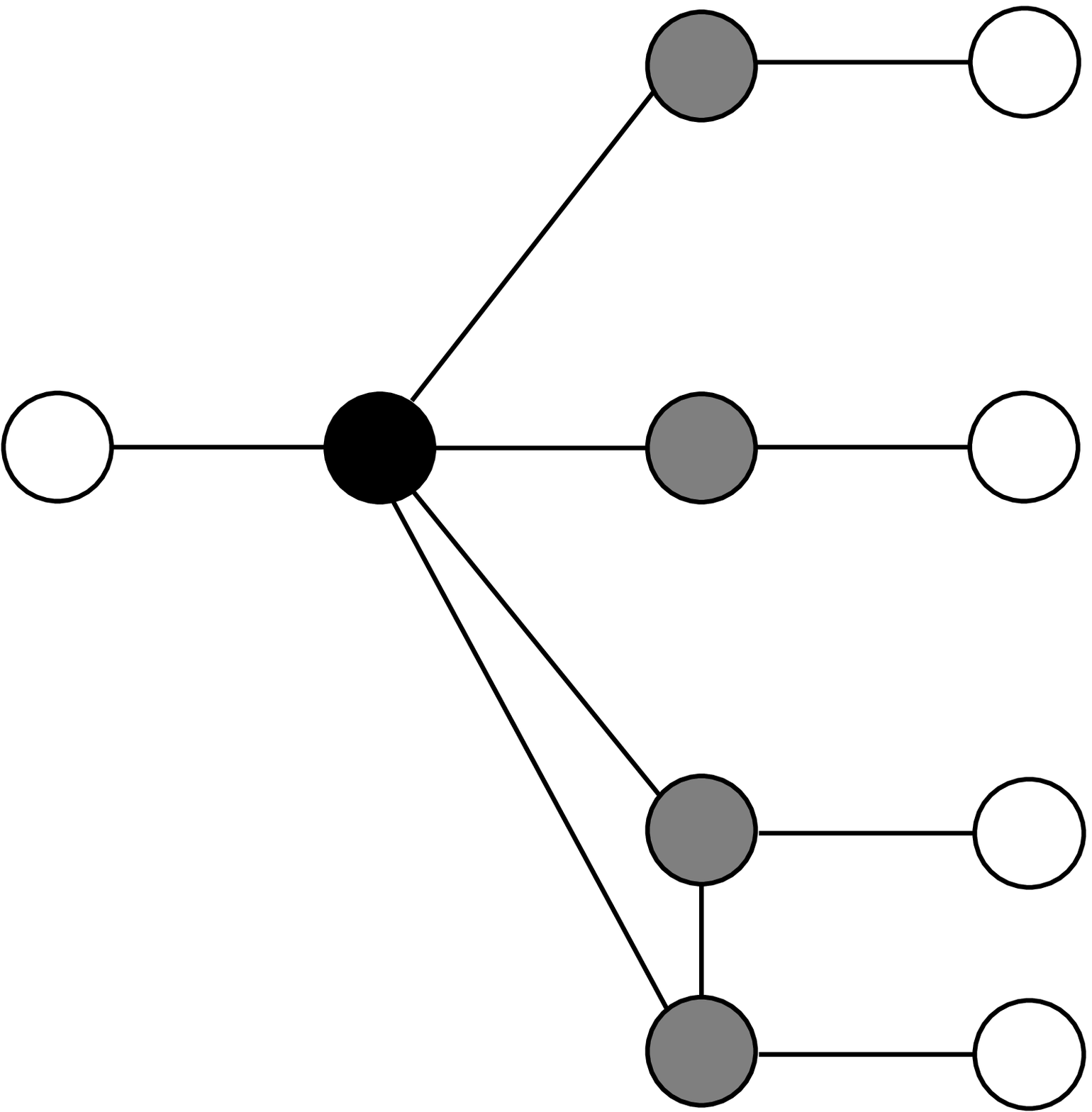, width=0.32\textwidth}
\put(-68,157){\footnotesize\sf Algeria}
\put(-71,104){\footnotesize\sf S.\hspace{-0.5mm} Africa}
\put(-60,52){\footnotesize\sf India}
\put(-24,157){\footnotesize\sf Tunisia}
\put(-21,104){\footnotesize\sf Israel}
\put(-22,52){\footnotesize\sf Oman}
\put(-153,106){\footnotesize\sf Pakistan}
\put(-77,-10){\footnotesize\sf Australasia}
\put(-15,-10){\footnotesize\sf Fiji}
\put(-111,106){\footnotesize\sf Core}
}
\hspace{50mm}
\parbox{0.5in}{\mbox{}\vspace{-1.5mm}\\
\epsfig{figure=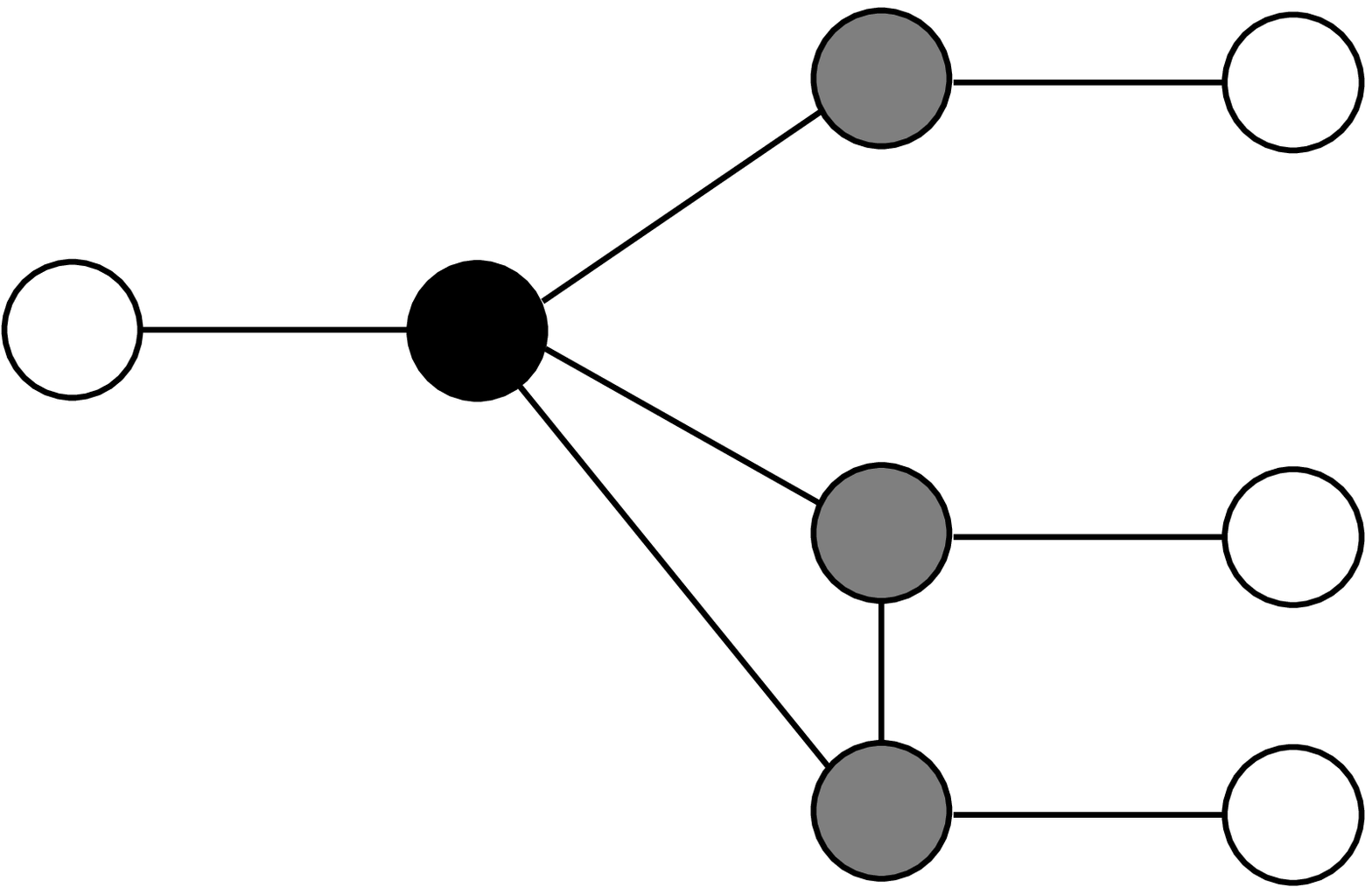, width=0.32\textwidth}
\put(-71,101){\footnotesize\sf S.\hspace{-0.5mm} Africa}
\put(-21,101){\footnotesize\sf Israel}
\put(-63,51){\footnotesize\sf India}
\put(-22,51){\footnotesize\sf Oman}
\put(-108,75){\footnotesize\sf Core}
\put(-153,75){\footnotesize\sf Pakistan}
\put(-77,-11){\footnotesize\sf Australasia}
\put(-15,-11){\footnotesize\sf Fiji}
}
\\\vspace{-2mm}\\
\caption{Clustering with threshold 125M USD: (a) CSP network; (b) Reduction of false twin vertices; (c) CSP structure.}
\end{figure}


As indicated earlier, the clustering criterion above already paves the way to illustrate 
some relation patterns; in a deeper analysis, however, it displays a severe limitation. Clustering
countries according to their amount of trade works well for (eventually defined) core clusters, and
also for some semiperipheries. But it does not accommodate the identification of
semiperiphery countries which, not trading a significant amount between themselves, display
however a similar (or even identical) connection pattern to the rest of the network. To incorporate
this, the criterion above should be combined with a similarity (or dissimilarity) measure
identifying countries with similar relation patterns. 

To illustrate this idea we first raise the trade threshold above to
500M USD. 
This yields a smaller cluster defined by the five East Asian countries
(trading more than 4.6 billion USD among themselves). 
Second, since we are dealing with a weighted network 
we define a dissimilarity criterion as follows: for each country we label each one
of its incident
edges with the percentage of trade that it carries,
computed over the country's total amount of trade.
This percentage is zero for absent edges, that is, 
for pairs of countries not
adjacent to each other. Denoting this percentage by $w_{ij}$ for 
the edge connecting vertices $i$ and $j$, 
the dissimilarity measure for countries $i$, $j$
is then defined as
$$\delta_{ij}=\sum_{i \neq k \neq j} |w_{ik}-w_{jk}|.$$
This means that two countries which have exactly the same connection
pattern to the rest of the network have a dissimilarity measure close
to zero (not exactly zero, in most cases, because even if the connections
are the same the percentages will typically
be different); on the contrary, if $i$ and $j$ are not adjacent and
do not have any neighbor in common
then the dissimilarity measure reaches 
the maximum value $\delta_{ij}=2$. 

Ignoring peripheries, we may now define new clusters (that
is, besides the
main one above) in terms of this dissimilarity measure: for instance, 
we may join together a set of countries into a single cluster
if the dissimilarities of all pairs within this set do not reach
a threshold of 1.0. Two non-trivial
clusters arise this way: the four Southeast Asian 
countries are joined into a single cluster (the six dissimilarities
range from 0.33 (Malaysia-Singapore) to 0.95 (Singapore-Philippines); the total internal
trade in this cluster reaches 585M USD),
and so do Australia and New Zealand (with a dissimilarity of 0.59; the trade among
themselves is 168M USD). The
remaining countries remain isolated. Note that none of these countries
reach, in any connection, the threshold of peer-to-peer trade of 500M USD defined above.

The quotient graph which results from this 
new clustering is displayed in Figure 7(a); now 
Sri Lanka is not adjacent to the core but
to the Southeast Asian cluster, via Singapore. As already depicted in this figure,
the five East Asian countries qualify again as a core, 
whereas the other clusters do not because of the eccentricity criterion.
The reductions of false twin vertices and of F-twin pairs yielding a CSP structure
can be found in Figure 7(b)-(c).
\begin{figure} \label{fig-AAOceania3}
\vspace{6mm}
\hspace{-10mm}
\parbox{0.5in}{
\epsfig{figure=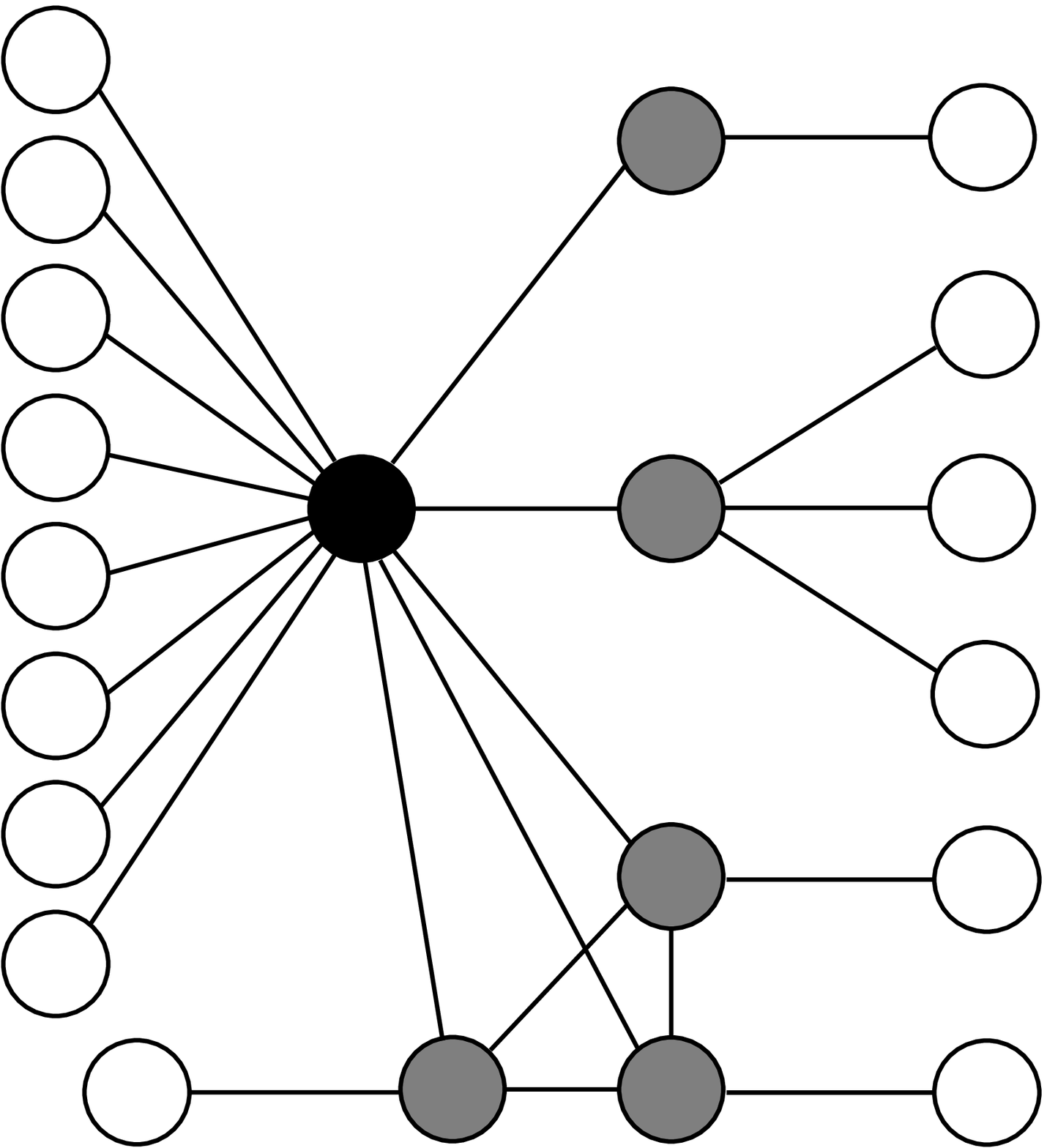, width=0.32\textwidth}
\put(-68,158){\footnotesize\sf Algeria}
\put(-74,105){\footnotesize\sf S.\hspace{-0.5mm} Africa}
\put(-56,49){\footnotesize\sf India}
\put(-64,-11){\footnotesize\sf Australasia}
\put(-111,-11){\footnotesize\sf Southeast}
\put(-97,-22){\footnotesize\sf Asia}
}
\hspace{50mm} 
\parbox{0.5in}{
\epsfig{figure=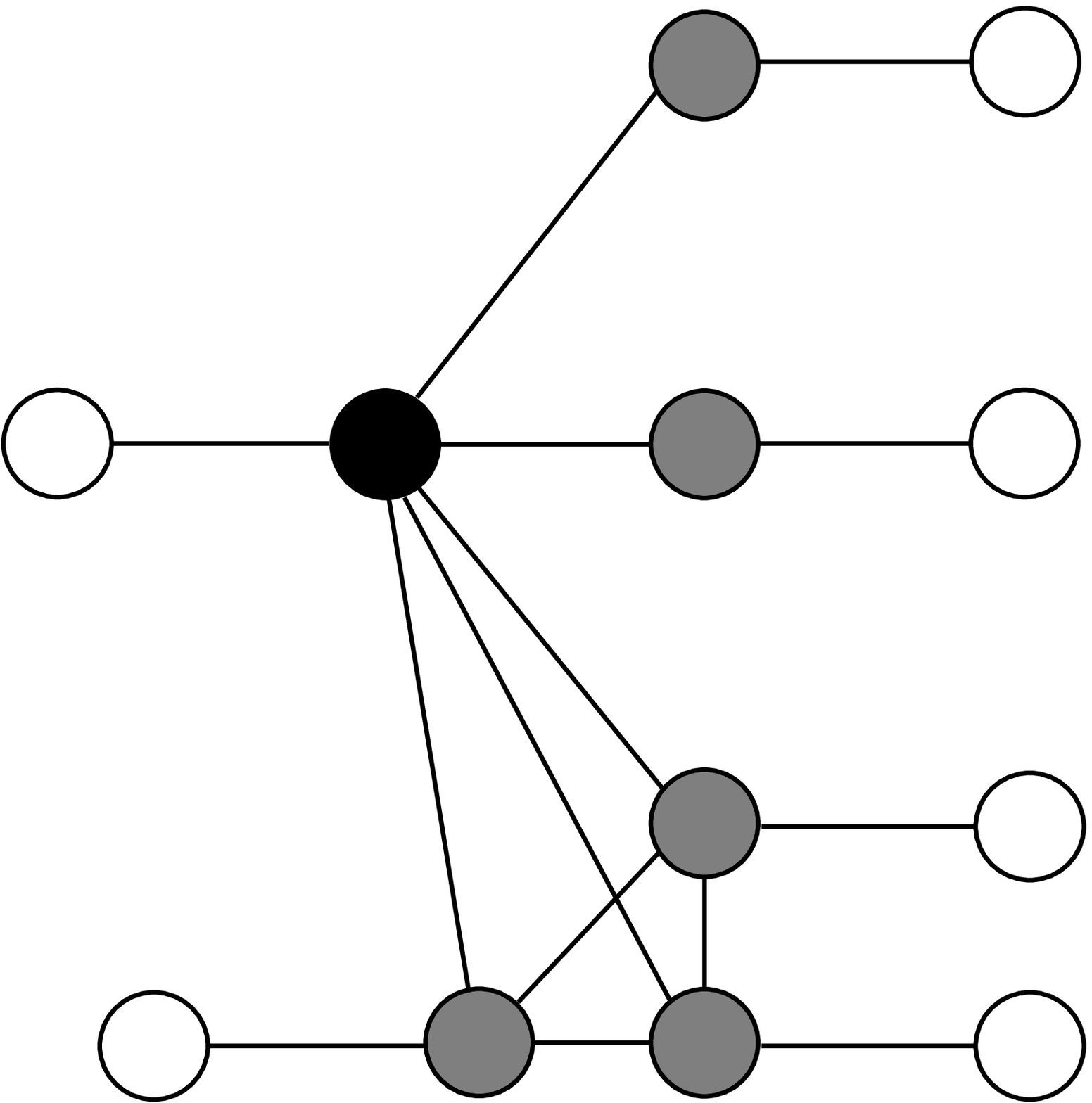, width=0.32\textwidth}
\put(-68,156){\footnotesize\sf Algeria}
\put(-71,102){\footnotesize\sf S.\hspace{-0.5mm} Africa}
\put(-56,50){\footnotesize\sf India}
\put(-23,156){\footnotesize\sf Tunisia}
\put(-20,102){\footnotesize\sf Israel}
\put(-21,50){\footnotesize\sf Oman}
\put(-153,102){\footnotesize\sf Pakistan}
\put(-64,-11){\footnotesize\sf Australasia}
\put(-14,-11){\footnotesize\sf Fiji}
\put(-110,-11){\footnotesize\sf Southeast}
\put(-96,-22){\footnotesize\sf Asia}
\put(-135,-11){\footnotesize\sf Sri}
\put(-141,-22){\footnotesize\sf  Lanka} 
\put(-111,102){\footnotesize\sf Core}
}
\hspace{50mm}
\parbox{0.5in}{\mbox{}\vspace{1mm}\\ 
\epsfig{figure=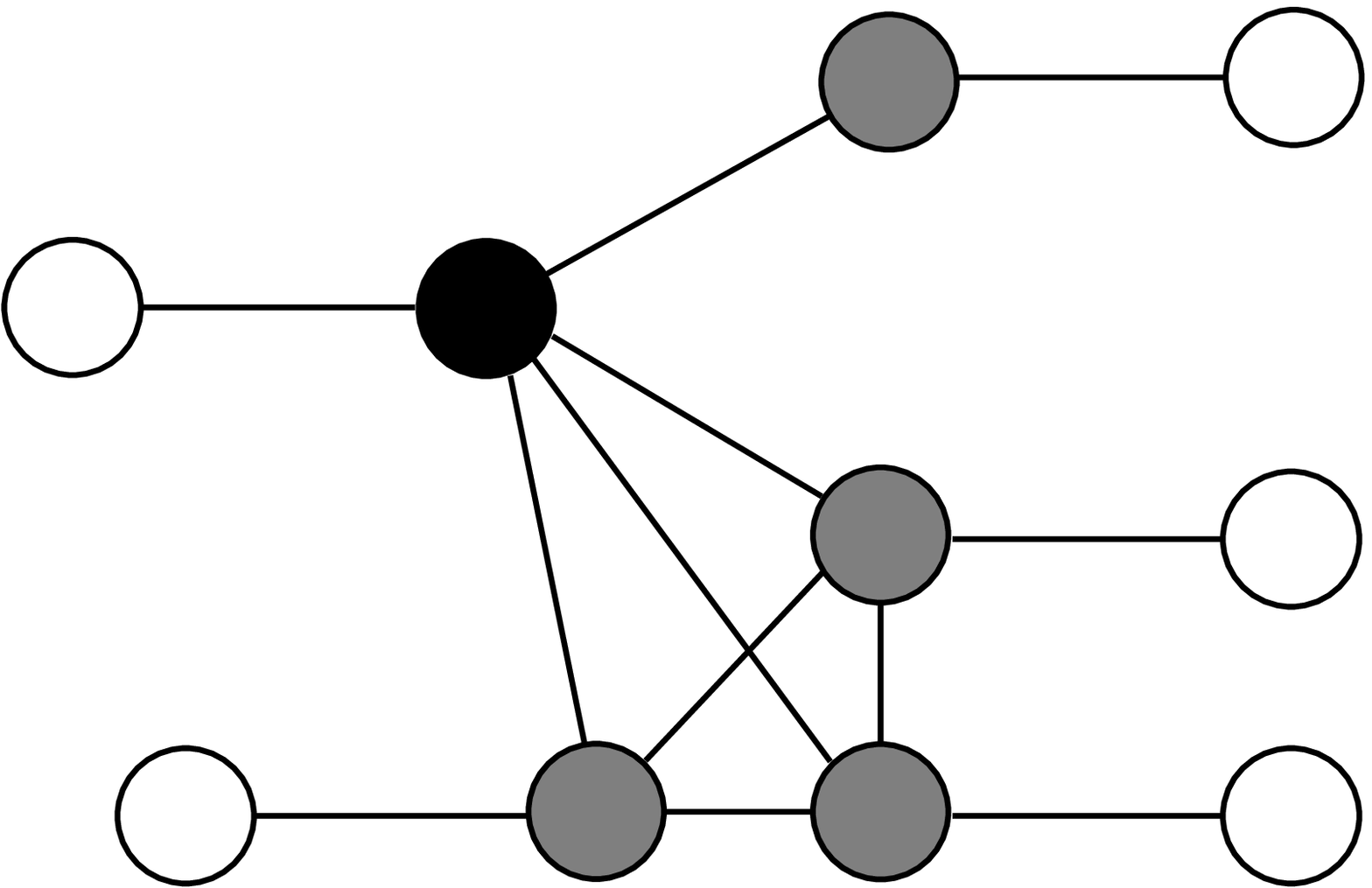, width=0.32\textwidth}
\put(-72,101.5){\footnotesize\sf S.\hspace{-0.5mm} Africa}
\put(-20,101.5){\footnotesize\sf Israel}
\put(-61,50){\footnotesize\sf India}
\put(-21,50){\footnotesize\sf Oman}
\put(-105,76){\footnotesize\sf Core}
\put(-153,76){\footnotesize\sf Pakistan}
\put(-64,-11){\footnotesize\sf Australasia}
\put(-14,-11){\footnotesize\sf Fiji}
\put(-110,-11){\footnotesize\sf Southeast}
\put(-96,-22){\footnotesize\sf Asia}
\put(-135,-11){\footnotesize\sf Sri}
\put(-141,-22){\footnotesize\sf  Lanka} 
}
\\\vspace{-2mm}\\
\caption{Adding a dissimilarity measure: (a) CSP network; (b) Reduction of false twin vertices; (c) CSP structure.}
\end{figure}

Finally, in order to further illustrate the eventual presence of other F-twin substructures, let us
ignore in Figure 7(c) 
the edge connecting India and Australasia: among 
the three semiperipheries at the bottom of this figure, this is clearly the one carrying less
trade (20,2M USD, whereas Southeast Asia trades 47,9M with India and over 177M with Australasia). The
resulting network is depicted in Figure 8(a).
Note that now the Australasia-Fiji and
India-Oman pairs become F-twins; they are isomorphic, disjoint and non-adjacent, and the connection
pattern to the remainder of the network is the same (both Australasia and India are connected to the
core and to Southeast Asia). We can therefore reduce this new relation pattern and the resulting
structure is displayed in  Figure 8(b). 
Worth clarifying is that the Australasia-Fiji
pair now
stands as the representative of this pattern, which is also met by the India-Oman pair.

\begin{figure} \label{fig-AAOceania4}
\vspace{6mm}
\hspace{15mm} 
\parbox{0.5in}{
\epsfig{figure=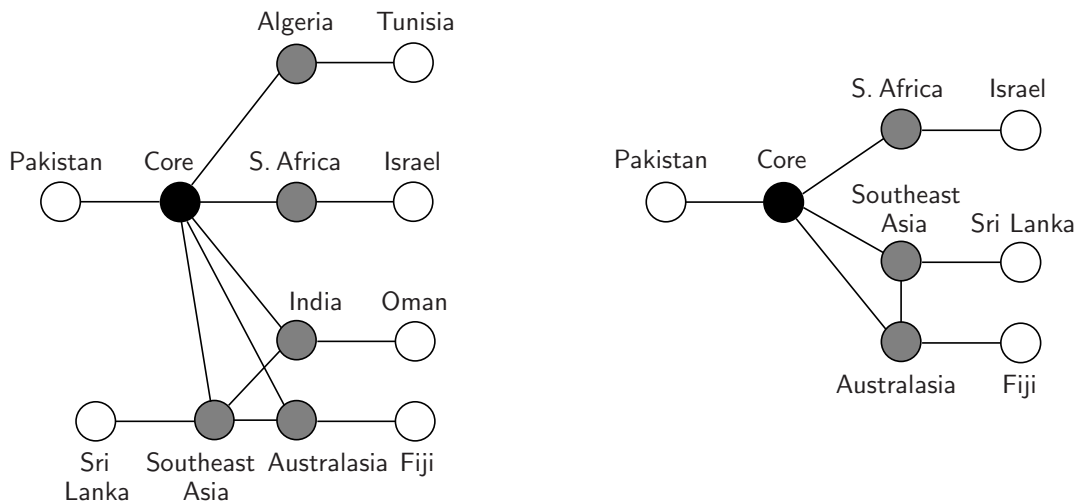, width=0.32\textwidth}
\put(-68,156){\footnotesize\sf Algeria}
\put(-71,103){\footnotesize\sf S.\hspace{-0.5mm} Africa}
\put(-56,50){\footnotesize\sf India}
\put(-23,156){\footnotesize\sf Tunisia}
\put(-20,103){\footnotesize\sf Israel}
\put(-21,50){\footnotesize\sf Oman}
\put(-162,103){\footnotesize\sf Pakistan}
\put(-64,-11){\footnotesize\sf Australasia}
\put(-14,-11){\footnotesize\sf Fiji}
\put(-110,-11){\footnotesize\sf Southeast}
\put(-96,-22){\footnotesize\sf Asia}
\put(-135,-11){\footnotesize\sf Sri}
\put(-141,-22){\footnotesize\sf  Lanka} 
\put(-111,103){\footnotesize\sf Core}
}
\hspace{65mm}
\parbox{0.5in}{\mbox{}\vspace{-5mm}\\ 
\epsfig{figure=AAOceania_secondcore_withoutFTwins.eps, width=0.32\textwidth}
\put(-72,101){\footnotesize\sf S.\hspace{-0.5mm} Africa}
\put(-20,101){\footnotesize\sf Israel}
\put(-72,61){\footnotesize\sf Southeast}
\put(-61,50){\footnotesize\sf Asia}
\put(-27,50){\footnotesize\sf Sri Lanka}
\put(-108,74){\footnotesize\sf Core}
\put(-162,74){\footnotesize\sf Pakistan}
\put(-78,-11){\footnotesize\sf Australasia}
\put(-15,-11){\footnotesize\sf Fiji}
}
\\\vspace{-2mm}\\
\caption{(a) The removal of the Australasia-India edge yields a new pair of F-twin subgraphs;
(b) Resulting CSP structure.}
\end{figure}

As indicated earlier in this section, the network here analyzed is intended
to illustrate the lines along which the results presented in this paper
can be applied to real problems. Future study should provide
a systematic analysis of clustering criteria in this context;
these criteria should combine density
and similarity measures. In a second step, quality measures defining
the extent to which the nodes in the quotient (clustered) graph  
may be classified either as cores, semiperipheries or peripheries
would indicate to what degree the network fits a
CSP structure. When a CSP structure is actually met, 
the twin notions here introduced make it possible to reduce identical
substructures, capturing the relation patterns depicted in the
network.

The example here considered suggests that the roadmap above is a promising
one. Note that the threshold
parameters within the aforementioned clustering criteria (involving e.g.\
the amount of trade or the degree of dissimilarity between countries)
has allowed for a 
progressive refinement of the clusters, providing gradually more detailed
information about the network structure. Indeed, the (say) giant core
in Figure 5(b) yields two clusters in Figure 6(b), namely
East/Southeast Asia 
and Australasia; in turn, the East-Southeast Asian core
is split in two in Figure 7(b). Accordingly,
the corresponding CSP structures in Figures 5(c), 6(c) and 7(c)
(with four, eight and ten nodes, respectively)
gradually display more detailed information about the
network structure. The network example here considered also shows
how different twin structures may be identified and reduced.
These include not only twin vertices but different
semiperiphery-periphery patterns:
compare e.g.\ in Figure 8(a) the 
Algeria-Tunisia and South Africa-Israel pairs, on the one hand,
and India-Oman and Australasia-Fiji, on the other.
Naturally, more complicated \spe patterns would arise in larger networks.

\section{Concluding remarks}
\label{sec-con}

Many problems related to twin subgraphs and to \csp structures remain
open for future study. We compile here some of them. First, 
the T-twin and F-twin notions for subgraphs introduced
in Sections \ref{sec-falsetwins} and \ref{sec-truetwins} have for sure 
a connection to automorphic and orbital equivalences, 
much as twin vertices arise in situations in which a 
transposition yields a graph automorphism. Note in this regard that, 
for vertices, the true
and false twin notions accommodate all possible cases of structurally
equivalent vertices, but for higher order subgraphs other
twin notions besides T-twins and F-twins
might be considered (for this reason we avoid
using the ``true'' and ``false'' labels for our T-twin and F-twin notions,
since the former labels seem to cover exhaustively all possible cases).
The classification of twin structures partially addressed in subsection
\ref{subsec-classification} also 
seems to have several potential extensions,
in particular connected to the interrelations between the
classification of different families of twin subgraphs.

Concerning the results considered in Section \ref{sec-csp},
it would be interesting to examine systematically to what extent
the set of actors (countries, companies, etc.) in
real social or economic networks can be clustered in a way that
matches some of the structures
displayed in subsection \ref{subsec-examples} after a suitable
reduction of twin patterns: the example discussed in Section
\ref{sec-ex} suggests a plan for future research in this
direction.
Motivated by the enumeration of CSP structures (cf.\ subsection
\ref{subsec-enumeration}),
several enumeration problems arise in connection to the
absence of twin substructures in graphs: specifically,
it would of interest to get a general enumeration formula 
for graphs without true twin vertices 
(or equivalently, in light of Theorem \ref{th-duality}, for graphs
without false
twin vertices), and also for graphs without any kind of T-twin
(or, analogously, F-twin) subgraphs. Closely related are the
problems of enumerating graphs without any kind of twin 
vertices, or without any kind of twin subgraphs. It also seems to be
worth studying other (say, layered) structures emanating
from greater parameter values in Definition \ref{defin-csp}, that
is, accommodating core eccentricities greater than two and/or periphery
degrees greater than one. All these topics are in the scope of future research.

\newpage

\noindent {\large\bf Appendix: The Asia-Africa-Oceania metal manufactures network}
\begin{center}
Table 3: Asia-Africa-Oceania metal manufactures trade in 1994 (from \cite{deNooy})\\
\mbox{} \vspace{0mm}\\
\begin{tabular}{ccc} \toprule
\ \ \ \ \ \ \ \ Country \#1 \ \ \ \ \ \ \ \ & \ \ \ \ \ \ Country \#2 \ \ \ \ \ \ & Trade  (thoushands of USD) \\ \midrule
China   &   Hong Kong   &   1482824\\
Japan   &   Thailand   &   894820\\
Japan   &   Korea   &   880295\\
China   &   Japan   &   630342\\
Malaysia   &   Singapore   &   484350\\
Japan   &   Malaysia   &   453463\\
Japan   &   Singapore   &   380454\\
Hong Kong   &   Japan   &   351919\\
Indonesia   &   Japan   &   200451\\
China   &   Korea   &   181392\\
Australia   &   New Zealand   &   168680\\
Japan   &   Philippines   &   138348\\
China   &   Singapore   &   135616\\
Japan   &   Australia   &   115283\\
Hong Kong   &   Singapore   &   110574\\
Singapore   &   Thailand   &   107720\\
China   &   Australia   &   90620\\
Australia   &   Indonesia   &   72387\\
Korea   &   Hong Kong   &   65315\\
Australia   &   Singapore   &   62392\\
Korea   &   Thailand   &   56160\\
Korea   &   Singapore   &   50098\\
Korea   &   Australia   &   45517\\
China   &   Thailand   &   44387\\
Australia   &   Malaysia   &   43068\\
Korea   &   Indonesia   &   41827\\
Indonesia   &   Malaysia   &   40291\\
China   &   Malaysia   &   39617\\
Singapore   &   Indonesia   &   39206\\
Malaysia   &   Thailand   &   37963\\
China   &   Indonesia   &   32817\\
India   &   Singapore   &   32130\\
Korea   &   Malaysia   &   31255\\
Japan   &   India   &   27655\\
Japan   &   South Africa   &   24555\\
\midrule 
\end{tabular}
\end{center}

\mbox{}\vspace{-18mm}\\

\flushright {\em (continued on next page)}\hspace{12mm}
\newpage
\begin{center}
\begin{tabular}{ccc} \toprule
\ \ \ \ \ \ \ \ Country \#1 \ \ \ \ \ \ \ \ & \ \ \ \ \ \ Country \#2 \ \ \ \ \ \ & Trade  (thoushands of USD) \\ \midrule
Hong Kong   &   Malaysia   &   24159\\
Hong Kong   &   Thailand   &   23642\\
Hong Kong   &   Philippines   &   23396\\
China   &   South Africa   &   23166\\
Singapore   &   Philippines   &   21744\\
Hong Kong   &   South Africa   &   21277\\
Australia   &   India   &   20366\\
China   &   Philippines   &   19865\\
Israel   &   South Africa   &   19183\\
Korea   &   South Africa   &   17826\\
Korea   &   Philippines   &   17031\\
India   &   Malaysia   &   15817\\
Japan   &   New Zealand   &   15470\\
Korea   &   Pakistan   &   15469\\
Thailand   &   Australia   &   14377\\
China   &   Egypt   &   14342\\
China   &   Pakistan   &   13953\\
Hong Kong   &   Australia   &   13644\\
China   &   New Zealand   &   12810\\
Hong Kong   &   Indonesia   &   12604\\
Singapore   &   Sri Lanka   &   12253\\
China   &   Algeria   &   11709\\
Australia   &   Fiji   &   10589\\
Japan   &   Pakistan   &   10388\\
China   &   Kuwait   &   9232\\
China   &   Jordan   &   8014\\
China   &   Morocco   &   7077\\
South Africa   &   Mauritius   &   6805\\
Algeria   &   Tunisia   &   6283\\
China   &   Bangladesh   &   5217\\
India   &   Oman   &   4151\\
Thailand   &   Seychelles   &   3179\\
South Africa   &   Reunion   &   2566\\
Japan   &   Madagascar   &   2042\\
\midrule 
\bottomrule
\end{tabular}
\end{center}

\end{document}